\documentclass[11pt,letterpaper]{amsart}
\usepackage{preamble}
\usepackage[createShortEnv]{proof-at-the-end}
\counterwithin{equation}{section}

\newcommand*{\vcenteredhbox}[1]{\begingroup
\setbox0=\hbox{#1}\parbox{\wd0}{\box0}\endgroup}

\begin{document}
\title[Vertex Connectivity of Friends-and-Strangers Graphs]{Vertex Connectivity of Friends-and-Strangers Graphs}
\author[Neil Krishnan]{Neil Krishnan}
\address[]{Massachusetts Institute of Technology, Cambridge, MA 02139}
\email{neilk301@mit.edu}
\author[Rupert Li]{Rupert Li}
\address[]{Department of Mathematics, Stanford University, Stanford, CA 94305}
\email{rupertli@stanford.edu}
\date{\today}

\begin{abstract}
Friends-and-strangers graphs, coined by Defant and Kravitz, are graphs on $n!$ vertices denoted by $\mathsf{FS}(X,Y)$ where $X$ and $Y$ are both graphs on $n$ vertices. The graph $X$ represents positions and edges denote adjacent positions while the graph $Y$ represents people and edges denote friendships. The vertex set of $\mathsf{FS}(X,Y)$ consists of all one-to-one placements of people on positions, and there is an edge between any two placements if it is possible to swap two people who are friends and on adjacent positions to get from one placement to the other. Previous papers have studied when $\mathsf{FS}(X,Y)$ is connected. In this paper, we consider when $\mathsf{FS}(X,Y)$ is $k$-connected where a graph is $k$-connected if it remains connected after removing any $k-1$ or less vertices. The connectivity of a graph is the largest $k$ for which it is $k$-connected.
We first consider $\mathsf{FS}(X,Y)$ when $Y$ is a complete graph or star graph.
We find tight bounds on their connectivity, proving their connectivity equals their minimum degree. We show that asymptotically similar conditions as the conditions mentioned by Bangachev are sufficient for $\mathsf{FS}(X,Y)$ to be $k$-connected. Finally, we consider when $X$ and $Y$ are independent Erd\H os--R\'enyi random graphs on $n$ vertices and edge probability $p_1$ and $p_2,$ respectively. We show that for $p_0 = n^{-1/2+o(1)},$ if $p_1p_2\geq p_0^2$ and $p_1,$ $p_2 \geq w(n) p_0$ where $w(n) \rightarrow 0$ as $n \rightarrow \infty,$ then $\mathsf{FS}(X,Y)$ is $k$-connected with high probability. This is asymptotically tight as we show that below an asymptotically similar threshold $p_0'=n^{-1/2+o(1)}$, the graph $\mathsf{FS}(X,Y)$ is disconnected with high probability if $p_1p_2 \leq (p_0')^2$.
\end{abstract}

\maketitle
\section{Introduction}
\subsection{Background}
The friends-and-strangers graph was introduced by Defant and Kravitz \cite{defant2021friends}.
\begin{definition}
    Given two simple $n$-vertex graphs $X=(V(X),E(X))$ and $Y=(V(Y),E(Y))$, define the \textit{friends-and-strangers graph} $\FS(X,Y)$ as follows.
    Its vertex set is the set of bijections $\sigma:V(X)\to V(Y)$, and $(\sigma,\tau)$ is an edge if and only if $\sigma=\tau\circ(i',j')$ for some $i',j'\in V(X)$ with $(i',j')\in E(X)$ and $(\sigma(i'),\sigma(j'))\in E(Y)$.
    In this case, we say $\sigma$ and $\tau$ differ by an \textit{$(X,Y)$-friendly swap}.
\end{definition}

We make the assumption that $V(X) = V(Y) = [n]$ by relabeling the vertices of $X$ and $Y$.
As a result, the bijections of $\FS(X,Y)$ can be thought of as permutations in the symmetric group $\mathfrak{S}_n.$
Also note that $\FS(X,Y)\cong\FS(Y,X)$ given by the isomorphism $f: V(\FS(X,Y)) \rightarrow V(\FS(Y,X))$ where $f(\sigma) = \sigma^{-1}.$
As a convention in this paper, we will generally mark vertices with primes to denote that they are in the $X$ graph and not mark them with primes to denote that they are in the $Y$ graph.

This setup for a friends-and-strangers graph has the following natural interpretation. Imagine there are $n$ people where people who are friends are connected in the graph $Y$ and $n$ positions where adjacent positions are connected in the graph $X$. Place the people on the positions with one person standing on each position. An $(X,Y)$-friendly swap is then a swap of two people such that they are friends and they are standing on adjacent positions. The graph $\FS(X,Y)$ then has a vertex for each possible placement with an edge connecting two placements if they differ by a single $(X,Y)$-friendly swap.

\begin{example}
Let $X = Y =P_3$ be the path graph on $3$ vertices. Then $\FS(X,Y)$ is the graph shown in \cref{fig:FSex}. This example was also mentioned in Defant and Kravitz \cite[Example 1.1]{defant2021friends}.
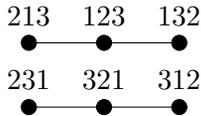
\begin{figure}[htbp]
    \centering
    \begin{tikzpicture}
        \draw (0,0) node[above,yshift=0.1 cm]{$213$}--(1,0) node[above,yshift=0.1 cm]{$123$}--(2,0) node[above,yshift=0.1 cm]{$132$};
        \filldraw (0,0) circle(0.1 cm) (1,0) circle(0.1 cm) (2,0) circle(0.1 cm);
    \end{tikzpicture}
    
    \vspace{0.1 cm}
    
    \begin{tikzpicture}
        \draw (0,0) node[above,yshift=0.1 cm]{$231$}--(1,0) node[above,yshift=0.1 cm]{$321$}--(2,0) node[above,yshift=0.1 cm]{$312$};
        \filldraw (0,0) circle(0.1 cm) (1,0) circle(0.1 cm) (2,0) circle(0.1 cm);
    \end{tikzpicture}

    \caption{The graph $\FS(P_3,P_3)$. The labels on the vertices mark the bijection from the positions to the people the vertex corresponds to, e.g., $312$ represents person $3$ is at position $1,$ person $1$ is at position $2,$ and person $2$ is at position $3.$}
    \label{fig:FSex}
\end{figure}
\end{example}
These graphs can model other combinatorial problems.
For example, as noted by Defant and Kravitz \cite{defant2021friends}, friends-and-strangers graphs generalize slide puzzles such as the famous 15-puzzle, which consists of the numbers $1$ through $15$, inclusive, placed on a $4\times4$ grid, leaving one cell empty. Then any number adjacent to the empty cell is allowed to move to the empty cell, leaving its original position empty, and the goal is to reach a predetermined configuration of numbers, canonically the numbers $1$ through $15$ in reading order, with the empty cell at the bottom right.
Such sliding puzzles can be represented by $\FS(X,Y)$ where $Y$ is a star graph with the empty cell as the central vertex adjacent to all other vertices, and $X$ being the adjacency graph of the sliding puzzle, i.e., $\operatorname{Grid}_{4\times 4}$ for the $15$-puzzle where $\operatorname{Grid}_{4\times 4}$ is $4\times 4$ grid graph.

In their initial paper, Defant and Kravitz \cite{defant2021friends} studied the basic properties of $\FS(X,Y)$ as well as the structures of $\FS(P_n,Y)$ and $\FS(C_n,Y)$ where $C_n$ is the cycle graph on $n$ vertices. 
In a later paper, Alon, Defant, and Kravitz \cite{alon2021extremal} considered the case when $X$ and $Y$ are both Erd\H os--R\'enyi random graphs $\mathcal{G}(n,p)$, identifying a threshold probability $p_{\text{gen}}=n^{-1/2+o(1)}$ above which $\FS(X,Y)$ is connected with high probability, and below which $\FS(X,Y)$ is disconnected with high probability.
Wang and Chen \cite{wang2023connectivity} considered the asymmetric variation of the problem where $X = \mathcal{G}(n,p_1)$ and $Y = \mathcal{G}(n,p_2)$.
Milojevi\'{c} \cite{milojevic2022connectivity} continued explorations of probabilistic questions, refining bounds on $p_{\text{gen}}$ and examining a generalization of friends-and-strangers graphs where $X$ and $Y$ have different numbers of vertices. Jeong \cite{jeong2022diameters} examined the diameters of friends-and-strangers graphs.

In addition, some extremal questions about friends-and-strangers graphs have been examined.
Alon, Defant, and Kravitz \cite{alon2021extremal} first posed the question of the smallest value of $d_n$ such that if $\delta(X),$ $\delta(Y) > d_n$ where $\delta(G)$ is the minimum degree of $G$, then $\FS(X,Y)$ is connected.
They proved that $3n/5-2 \leq d_n$ for all $n$ and $d_n \leq 9n/14+2$ when $n \geq 16.$ Bangachev \cite{bangachev2022asymmetric} refined the upper bound by proving that $d_n \leq \lceil 3n/5 \rceil,$ showing that $d_n = 3n/5 + O(1)$.

This paper generalizes some of these extremal and probabilistic results, applying it to vertex connectivity.
Recall a graph is \textit{$k$-connected} if we can remove any $k-1$ or less vertices and the remaining induced graph is still connected, or equivalently, for any two vertices $u$ and $v$, there exist $k$ paths between $u$ and $v$ whose vertices are pairwise disjoint (other than vertices $u$ and $v$ themselves). The \textit{connectivity} of a graph is the maximum $k$ for which it is $k$-connected.

One reason for studying vertex connectivity is it relates very naturally to sliding block puzzles like the 15-puzzle as a vertex connectivity of $k$ implies that for any $k-1$ forbidden configurations of the puzzle, there still exists a path between any two of the remaining configurations.

\subsection{Main Results}
In this section, we will establish the main results that will be proved in this paper. We start by turning the condition of vertex connectivity into an easier condition of the existence of disjoint paths between $\sigma$ and $\sigma \circ (u',v')$ for some $u',v' \in V(X).$ To do this, we find the connectivity of $\FS(X,K_n).$
\begin{theorem} \label{knconnectivity}
    For $n\geq 3$, if $G = \FS(X,K_n)$ where $X$ is connected, its connectivity is equal to its minimum degree.
\end{theorem}
Through this theorem and \cref{kdisjointpaths}, a result which shows how to combine disjoint paths between vertices, we have the following proposition which simplifies the condition of vertex connectivity.
\begin{proposition} \label{exchangeability}
    Let $n \geq k+1.$ Assume that for any $\sigma \in V(\FS(X,Y))$ and $u',v' \in V(X)$ where $(u',v') \in E(X),$ there exist $k$ disjoint paths between $\sigma$ and $\sigma \circ (u',v').$ Then, $\FS(X,Y)$ is $k$-connected.
\end{proposition}

We also extend Wilson's Theorem, which we formally state in \cref{section:preliminaries}, to vertex connectivity.
\begin{theorem} \label{genwilson}
    For $n\geq 3$, if $G = \FS(X,\Star_n)$ where $X$ is connected, the connectivity of its connected components is equal to its minimum degree. 
\end{theorem}
We use \cref{exchangeability,genwilson} in the proofs of the two theorems below which show two sets of sufficient conditions for $\FS(X,Y)$ to be $k$-connected.
\begin{theorem}\label{kban1}
    Let $X$ and $Y$ be two graphs on $n$ vertices satisfying
    \begin{itemize}
        \item $\delta(X), \delta(Y)>n/2,$
        \item $2\min\{\delta(X), \delta(Y)\}+3\max\{\delta(X),\delta(Y)\} \geq 3n+2k-4,$
    \end{itemize}
    where $k \geq 2.$ Then for sufficiently large $n$, namely $n \geq 5(1+(k-1)(k+6)+2k-3),$ we have $\FS(X,Y)$ is $k$-connected.
\end{theorem}
Notice that for fixed $k,$ the bounds are asymptotically the same as the bound conditions for classical $1$-connectivity proven by Bangachev \cite{bangachev2022asymmetric}; we discuss this and other previous results in more detail in \cref{section:preliminaries}.
Moreover, we know there do not exist $\alpha < 3$ for which the right hand side of the second bound condition becomes $\sim \alpha n$ as a result of Proposition 1.6 of \cite{bangachev2022asymmetric}. Conjecture 8.2 in \cite{bangachev2022asymmetric} conjectures the bounds of \cref{kban1} are tight in the sense that there is a counter example for minimum degrees not satisfying the bounds.
We also have the following result which eliminates the $\delta(X),$ $\delta(Y)>n/2$ condition at the cost of a worsened second bound condition.

\begin{theorem}\label{kban2}
    Let $X$ and $Y$ be two connected graphs on $n$ vertices satisfying
    \begin{itemize}
        \item $\delta(X) + \delta(Y) \geq n+k-1,$
        \item $\min\{\delta(X), \delta(Y)\}+2\max\{\delta(X),\delta(Y)\} \geq 2n,$
    \end{itemize}
    where $k \geq 2.$ Then $\FS(X,Y)$ is $k$-connected.
\end{theorem}
The bound conditions are asymptotically the same as the corresponding bounds for classical $1$-connectivity \cite{bangachev2022asymmetric}.

We also resolve the following question of Milojevi\' c with \cref{kconnectivitywhp}, showing that under the conditions of the question, for all $k$, for sufficiently large $n$, we have that $\FS(X,Y)$ is $k$-connected.
\begin{question} [\protect{\cite[Question 6.3]{milojevic2022connectivity}}]\label{milojevicquestion}
    Let $p(n) = n^{-1/2+o(1)}$ and let $X$ and $Y$ be random graphs in $\mathcal{G}(n,p).$ For which values of $k$ is the graph $\FS(X,Y)$ $k$-connected with high probability?
\end{question}
\begin{theorem}\label{kconnectivitywhp}
    Fix some small $\varepsilon > 0$ and a positive integer $k>1.$ Let $X$ and $Y$ be independently-chosen random graphs in $\mathcal{G}(n,p_1)$ and $\mathcal{G}(n,p_2)$ where $p_1=p_1(n)$ and $p_2=p_2(n)$ both depend on $n.$ If
    $$p_1p_2\leq \varepsilon\frac{\log n}{n},$$
    then $\FS(X,Y)$ is disconnected with high probability. If $p_1p_2 \geq p_0^2$ and $p_1,p_2\geq p_0/\ell$ where
    $$p_0= \frac{\exp((k+7)/4\cdot(\log n)^{2/3})}{n^{1/2}},$$
    and $\ell = \frac{(\log n)^{2/3}}{2},$ then $\FS(X,Y)$ is $k$-connected with high probability.
\end{theorem}
As a result for all $k,$ assuming the conditions of \cref{milojevicquestion}, $\FS(X,Y)$ is $k$-connected with high probability.

\subsection{Overview}
In \cref{section:preliminaries}, we provide the notation and terminology used throughout the paper and some related results in $1$-connectivity, i.e., classical connectivity. 
Then in \cref{section:exchangeability}, we show how to extend the idea of exchangeability \cite{alon2021extremal} in classical connectivity to vertex connectivity.
In \cref{section:wgen2}, we prove \cref{genwilson} which extends Wilson's theorem to determine when  $\FS(X,\Star_n)$ is $k$-connected.
In \cref{section:bangachev}, we prove \cref{kban1} and \cref{kban2}, and in \cref{section:random}, we prove \cref{kconnectivitywhp}.

\section{Preliminaries}\label{section:preliminaries}
\subsection{Notation and Terminology}
We assume that all graphs are simple if not otherwise stated. For a graph $G,$ let $V(G)$ denote the vertex set of $G,$ and let $E(G)$ denote the edge set of $G.$ The following graphs with vertex set $[n]$ will reoccur throughout the paper.
\begin{itemize}
    \item $K_n$ - the complete graph which has edge set $\{ (i,j) : 1 \leq i < j \leq n\};$
    \item $\Star_n$ - the star graph which has edge set $\{ (i,n) :1 \leq i < n\}.$
    \item $\Star^+_n$ - a graph with edge set $ \{ (i,n) : 1 \leq i \leq n-1 \} \cup \{(1,2)\}.$
    \item $P_n$ - the path graph which has edge set $\{ (i,i+1) :1 \leq i < n\}.$ Note that we use $n$ for the number of vertices instead of the number of edges.
    \item $C_n$ - the cycle graph which has edge set $\{ (i,i+1): 1 \leq i < n\}\cup \{(1,n)\}.$
\end{itemize}
We use additional terminology to describe a graph $G.$
\begin{itemize}
    \item For some subset $V_0$ of $V(G),$ the \textit{induced subgraph} denoted by $G\lvert_{V_0}$ is the graph $H$ where $V(H) = V_0$ and $E(H) = \{(i,j) : i,j \in V_0,\, (i,j) \in E(G)\}.$
    \item Define $\delta(G)$ to be the \textit{minimum degree} of $G.$
    \item We say $G$ is \textit{$k$-connected} for any set of $k-1$ or less vertices, after removing them, the remaining induced graph is connected. Specifically if $G$ satisfies this property for $k = 2,$ we refer to $G$ as \textit{biconnected}. Note that we recover the classical definition of a graph being connected with the case $k=1$.
    We define the \textit{connectivity} of $G$ to be the maximum $k$ for which $G$ is $k$-connected.
    \item A \textit{cut vertex} of $G$ is $v \in V(G)$ such that $G\lvert_{V(G)\backslash \{v\}}$ has more connected components than $G.$
    \item Two paths, $P_1$ and $P_2,$ in $G$ are \textit{disjoint} if $P_1$ and $P_2$ share no other vertices other than their starting vertex or ending vertex.
    \item For a graph $G,$ we define the \textit{open neighborhood} $N(v)$ for vertex $v \in V(G)$ to be the set of vertices adjacent to $v.$ We further define the \textit{closed neighborhood} $N[v] = N(v) \cup \{v\}.$
    \item For a set $S \subseteq V(G),$ we define $N(S)$ to be the set of vertices in $G$ that are not in $S$ but which are adjacent to a vertex of $S.$
\end{itemize}

One property of friends-and-strangers graphs that we will use implicitly throughout this paper is how $\FS(X,Y)$ changes on restricting to a subgraph of $X.$
\begin{proposition}[\protect{\cite[Proposition 2.1]{defant2021friends}}]
    Let $X, \widetilde{X},Y,\widetilde{Y}$ be graphs on $n$ vertices such that $\widetilde{X}$ is a subgraph of $X$ and $\widetilde{Y}$ is a subgraph of $Y$. Then $\FS(\widetilde{X},\widetilde{Y})$ is a subgraph of $\FS(X,Y).$
\end{proposition}
This proposition holds more generally as long as $\widetilde{X}$ and $\widetilde{Y}$ have the same number of vertices which is not necessarily the same as the number of vertices $X$ and $Y$ have.

\subsection{Prior Results on Connectivity of $\FS(X,Y)$}
We summarize relevant work on friends-and-strangers graphs.
\subsubsection{Challenges to Connectivity}
There are two major obstructions to connectivity in $\FS(X,Y).$ The first one is when there exists a cut vertex in both $X$ and $Y.$
\begin{proposition}[\protect{\cite[Proposition 2.6]{defant2021friends}}]\label{cutv}
    Let $X$ and $Y$ be connected graphs, each on $n\geq 3$ vertices. Suppose $x_0 \in X$ and $y_0 \in Y$ are cut vertices such that the connected components of $X\lvert_{V(x)\backslash \{x_0\}}$ are $X_1,\ldots, X_r$ and the connected components of $Y\lvert_{V(Y)\backslash\{y_0\}}$ are $Y_1,\ldots, Y_s.$ Let $\mathcal{M}=\mathcal{M}(X,Y,x_0,y_0)$ denote the set of $r\times s$ matrices with nonnegative integer entries in which the $i$th row sums to $\lvert V(X_i)\lvert$ and the $j$th column sums to $\lvert V(Y_j) \lvert.$ Then the number of connected components of $\FS(X,Y)$ is at least $\lvert \mathcal{M} \lvert.$
\end{proposition}
Let $X$ and $Y$ be graphs on $n\geq 3$ vertices such that both graphs have a cut vertex. As mentioned in the discussion following Proposition 2.6 in \cite{defant2021friends},because $\lvert\mathcal{M}\rvert  \geq 2$ since $r,s \geq 2,$ we know $\FS(X,Y)$ must be disconnected.

In addition, if $X$ and $Y$ are bipartite graphs, we run into a parity obstruction which forces $\FS(X,Y)$ to be disconnected.
For a permutation $\sigma$, let $\sgn(\sigma)$ be the sign of $\sigma$, which is 1 if $\sigma$ is even, i.e., a product of an even number of transpositions, and $-1$ if $\sigma$ is odd.
\begin{proposition}[\protect{\cite[Proposition 2.7]{defant2021friends}}] \label{bip}
    Let $X$ and $Y$ be bipartite graphs on $n\geq 3$ vertices with vertex bipartitions $V(X) = A_X \sqcup B_X$ and $V(Y) = A_Y \sqcup B_Y.$ Given a permutation $\sigma,$ let
    \[ p(\sigma) = \lvert \sigma(A_X) \cap A_Y \vert + \frac{\sgn(\sigma)+1}{2}. \]
    If the permutations $\sigma$ and $\tau$ are in the same connected component of $\FS(X,Y),$ then $p(\sigma)$ and $p(\tau)$ have the same parity.
\end{proposition}
Note that both parities can be achieved as mentioned in the last paragraph of \cite[Proposition 2.7]{defant2021friends}. Consider any permutation $\sigma.$ Notice that for $n\geq 3,$ at least one of $A_X$ or $B_X$ has more than $1$ vertex. Without loss of generality, let $\lvert A_X \rvert \geq 2.$ Let $i',j' \in A_X$ where $i' \neq j',$ and let $\sigma' = \sigma \circ (i',j').$ Notice that $\lvert \sigma(A_X) \cap A_Y \rvert = \lvert \sigma'(A_X) \cap A_Y \lvert$ since the swap from $\sigma$ to $\sigma'$ was contained within $A_X,$ but $\sgn(\sigma) \neq \sgn(\sigma')$ because the number of transpositions changed by $1.$ Therefore both parities can be achieved, meaning $\FS(X,Y)$ is disconnected.
\subsubsection{Star Graphs}
As a result of the previous propositions, if we consider the special case $X = \Star_n,$ we can completely classify the graphs $X$ such that $\FS(X,\Star_n)$ is connected through Wilson's Theorem stated below.
\begin{theorem}[\protect{\cite[Theorem 1]{wilson1974graph}}] \label{wilson}
    Suppose that $X$ is a graph on $n \geq3$ vertices satisfying the following properties:
    \begin{itemize}
        \item $X$ is biconnected,
        \item $X$ is not bipartite,
        \item $X$ is not isomorphic to $C_n$ for $n \geq 4$,
        \item $X$ is not isomorphic to the graph $\theta_0$ depicted in \cref{fig:theta0}.
        \begin{figure}[htbp]
            \centering
            \begin{tikzpicture}
                \draw (-1,0)--(0,0)--(1,0)--(0.5,0.866)--(-0.5,0.866)--(-1,0)--(-0.5,-0.866)--(0.5,-0.866)--(1,0);
                \filldraw (-1,0) circle(0.1 cm) (0,0) circle(0.1 cm) (1,0) circle(0.1 cm) (0.5,0.866) circle(0.1 cm) (-0.5,0.866) circle(0.1 cm) (-0.5,-0.866) circle(0.1 cm) (0.5,-0.866) circle(0.1 cm);
            \end{tikzpicture}
            \caption{The graph $\theta_0$.}
            \label{fig:theta0}
        \end{figure}
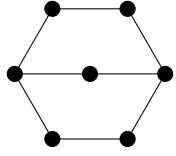
    \end{itemize}
    Then $\FS(X,\Star_n)$ is connected.
\end{theorem}
We will refer to graphs $Y$ satisfying all four above properties as being \textit{Wilsonian}. Because of the previous propositions, notice that if $Y$ doesn't satisfy the first two properties of Wilsonian graphs, for $n\geq 3,$ we have that $\FS(X,Y)$ must be disconnected. In addition, as Defant and Kravitz \cite{defant2021friends} mentioned before their Theorem 2.5, $\FS(\Star_n, C_n)$ must have $(n-2)!$ connected components, meaning that $Y$ not satisfying the third property are also disconnected.
Finally, $\FS(\Star_7,\theta_0)$ has $6$ connected components, showing that we have a complete characterization of connected $\FS(\Star_n,Y)$ graphs for $n \geq 3$.


\subsubsection{Minimum Degree Conditions}
Bangachev \cite{bangachev2022asymmetric} studied connectivity of $\FS(X,Y),$ using notions of $(X,Y)$-exchangeability developed in \cite{alon2021extremal}.
In particular, vertices $u,v \in V(Y)$ are \textit{$(X,Y)$-exchangeable} if for every permutation $\sigma,$ there exists a series of $(X,Y)$-friendly swaps which turns $\sigma$ into $(u,v)\circ \sigma.$
\begin{lemma} [\protect{\cite[Lemma 2.9]{alon2021extremal}}]\label{exchangeable}
    Let $X$ and $Y$ be two graphs on $n$ vertices. Assume that for all vertices $u,v \in V(Y),$ and all permutations, $\sigma,$ where $(\sigma^{-1}(u),\sigma^{-1}(v)) \in E(X),$ we have that $u$ and $v$ are $(X,Y)$-exchangeable. Then $\FS(X,Y)$ is connected.
\end{lemma}
This notion is powerful because it reduces the global problem of connectivity into a more localized problem of swapping just two values in a permutation.
Using these ideas as well as \cref{wilson}, Bangachev proved two asymmetric conditions on connectivity of $\FS(X,Y).$ 
\begin{theorem} [\protect{\cite[Theorem 1.4]{bangachev2022asymmetric}}]\label{ban1}
    Suppose that $X$ and $Y$ are two graphs on $n\geq 6$ vertices satisfying:
    \begin{itemize}
        \item $\delta(X), \delta(Y) > n/2,$
        \item $2\min\{\delta(X),\delta(Y)\}+3\max\{\delta(X),\delta(Y)\} \geq 3n.$
    \end{itemize}
    Then $\FS(X,Y)$ is connected.
\end{theorem}
\begin{theorem} [\protect{\cite[Theorem 1.5]{bangachev2022asymmetric}}] \label{ban2}
    Suppose that $X$ and $Y$ are two graphs on $n$ vertices satisfying:
    \begin{itemize}
        \item $X$ and $Y$ are both connected,
        \item $\min\{\delta(X),\delta(Y)\}+2\max\{\delta(X),\delta(Y)\} \geq 2n.$
    \end{itemize}
    Then $\FS(X,Y)$ is connected.
\end{theorem}
Notice that $d_n \leq \lceil 3n/5 \rceil$ follows from \cref{ban1}. \cref{kban1} and \cref{kban2} extend these results to $k>1.$

\subsubsection{Probability Thresholds} Alon, Defant, and Kravitz \cite{alon2021extremal} studied the case when $X$ and $Y$ are graphs independently chosen in $\mathcal{G}(n,p).$ They found a probability threshold below which $\FS(X,Y)$ is disconnected with high probability and above which $\FS(X,Y)$ is connected with high probability.
\begin{theorem}[\protect{\cite[Theorem 1.1]{alon2021extremal}}]\label{connectivitywhp}
    Fix some small $\varepsilon > 0$ and a positive integer $k>1.$ Let $X$ and $Y$ be independently-chosen random graphs in $\mathcal{G}(n,p),$ where $p=p(n)$ depends on $n.$ If
    $$p\leq \frac{2^{-1/2}-\varepsilon}{n^{1/2}},$$
    then $\FS(X,Y)$ is disconnected with high probability. If
    $$p\geq \frac{\exp(2(\log n)^{2/3})}{n^{1/2}},$$
    then $\FS(X,Y)$ is connected with high probability.
\end{theorem}
Milojevi\' c \cite{milojevic2022connectivity} improved the disconnected with high probability bound by showing the following theorem.
\begin{theorem} [\protect{\cite[Theorem 1.2]{milojevic2022connectivity}}] \label{milojevicdisconnectedwhp}
    There exists a constant $\varepsilon > 0$ with the following property. For a large positive integer $n$ and
    $$p<\varepsilon \left( \frac{\log n}{n} \right)^{1/2},$$
    if we choose random graphs $X$ and $Y$ in $\mathcal{G}(n,p)$ independently, $\FS(X,Y)$ has an isolated vertex with high probability.
\end{theorem}
Wang and Chen \cite{wang2023connectivity} generalized \cref{connectivitywhp} to an asymmetric variant.
\begin{theorem} [\protect{\cite[Theorem 1.3]{wang2023connectivity}}] \label{wangchenconnectedwhp}
Fix some small $\varepsilon > 0,$ and let $X$ and $Y$ be independently chosen random graphs in $\mathcal{G}(n,p_1)$ and $\mathcal{G}(n,p_2),$ respectively, where $p_1 = p_1(n)$ and $p_2=p_2(n)$ depend on $n$. Let $p_0 = \exp(2(\log n)^{2/3})n^{-1/2}.$ If either
$$p_1p_2 \leq \frac{(1-\varepsilon)/2}{n} \text{ and } p_1,p_2 \gg \frac{\log n}{n},$$
or
$$\min\{p_1,p_2\} \leq \frac{\log n + c(n)}{n} \text{ for some }c(n) \rightarrow -\infty,$$
where $\ell = \frac{(\log n)^{2/3}}{2},$ then $\FS(X,Y)$ is disconnected with high probability. If $$p_1p_2 \geq p_0^2 \quad \text{and} \quad p_1,p_2\geq \frac{1}{\ell}p_0,$$ then $\FS(X,Y)$ is connected with high probability.
\end{theorem}

An extension of \cref{milojevicdisconnectedwhp} implies the disconnected side of \cref{kconnectivitywhp}. The connected side on the other hand is an extension of \cref{wangchenconnectedwhp}.

\subsection{Menger's Theorem}
Menger's theorem is a connection between two different ideas of connectivity.
\begin{theorem}[\protect{\cite{menger1927theorem}}]
    Let $G$ be a graph and $A,B \subseteq V(G).$ The minimum number of vertices $v \in V(G) \setminus A\cup B$ that need to be removed to disconnect $A$ and $B$ is equal to the maximum number of pairwise disjoint paths between $A$ and $B.$
\end{theorem}
As a corollary, we see that when $G$ has more than $k$ vertices, $G$ is $k$-connected if and only if there exist $k$ disjoint paths between any two vertices of $G$ \cite[Theorem 3.3.6]{diestelblocks}.
When $G$ has $k$ vertices, notice that $G$ is $k$-connected only when $G$ is the complete graph on $k$ vertices.

Note that when we are discussing $\FS(X,Y)$ where $X$ and $Y$ are graphs on $n$ vertices, we will generally restrict to $n\geq 3.$ As a result notice that $\FS(X,Y)$ cannot be a complete graph since the number of vertices in $\FS(X,Y)$ is $n!$ but the maximum possible degree of a vertex in $\FS(X,Y)$ is $n(n-1)/2$ which occurs when $X=Y=K_n$ since then any of $n(n-1)/2$ transpositions can be made on a permutation. Because $n\geq 3,$ we know $n!>n(n-1)/2+1$ meaning $\FS(X,Y)$ is not complete.

Therefore, when $n\geq 3,$ we know $\FS(X,Y)$ is $k$-connected if and only if there exist $k$ disjoint paths between any two vertices of $\FS(X,Y).$ We use this fact implicitly throughout the paper.

\subsection{Atomic parts}
We now introduce atomic parts which we will also sometimes refer to as atoms. Let $G$ be a graph. Define a cutset as follows.
\begin{definition}
    A cutset, $C,$ is a set of vertices in $G$ such that $G\lvert_{V(G)\setminus C}$ is disconnected.
\end{definition}
We further define $\mc{C}$ to be the set of cutsets of cardinality $k,$ the connectivity of $G$, i.e., the minimum cardinality.
Let $\rho(G)$ be the size of the smallest connected component out of all components generated by removing some cutset in $\mc{C}$ from $G.$
\begin{definition}
    An \textit{atomic part} is a connected component, $A,$ generated by removing a cutset of $\mc{C}$ such that $\lvert V(A) \rvert = \rho(G).$
\end{definition}
From Watkins \cite{watkins1970atomicparts}, we have the following theorems about atomic parts which will be relied on heavily throughout the paper.
\begin{theorem} [\protect{\cite[Theorem 1]{watkins1970atomicparts}}] \label{disjointatomic parts}
    In a connected graph, atomic parts are disjoint.
\end{theorem}
\begin{lemma} [\protect{\cite[Lemma 3.5]{watkins1970atomicparts}}] \label{adjacentatomic parts}
    Let $A$ and $A'$ be two atomic parts. If there exists a vertex in $V(A')\cap N(V(A)),$ then $V(A') \subseteq N(V(A)).$
\end{lemma}

\section{Exchangeability}\label{section:exchangeability}
As mentioned, a key part of connectivity proofs is $(X,Y)$-exchangeability which converts connectivity into a local problem. This section extends those ideas to vertex connectivity.
\begin{proposition}\label{kdisjointpaths}
    Let $k \geq 2$ be an integer. Let $u$ and $v$ be two vertices in a graph $G.$ Let $a_i \geq 1$ be integers for $1 \leq i \leq k.$ Let $x_{i,j}$ for $1 \leq j \leq a_i$ and $1 \leq i \leq k$ be distinct vertices in $G.$ Let $x_{i,0} = u$ and $x_{i,a_i+1} = v$ for all $1\leq i \leq k.$ If there exist $k$ disjoint paths between $x_{i,j}$ and $x_{i,j+1}$ for $1 \leq i \leq k$ and $0 \leq j \leq a_i,$ then there exist $k$ disjoint paths between $u$ and $v.$
\end{proposition}

\begin{proof}
    Proving that there exist $k$ disjoint paths between $u$ and $v$ is equivalent to showing that removing any $k-1$ vertices which are not $u$ or $v$ from $G$ does not disconnect $u$ and $v.$ Notice that if we remove $k-1$ vertices from $G,$ there must exist some $1 \leq i \leq k$ such that for all $x_{i,j}$ for $1 \leq j \leq a_i$ have not been removed from $G.$ Because there are $k$ disjoint paths from $x_{i,j}$ to $x_{i,j+1}$ for $0 \leq j \leq a_i,$ we know that $x_{i,j}$ and $x_{i,j+1}$ could not have been disconnected by removing $k-1$ vertices. Therefore, there exists a path between $u$ and $v.$
\end{proof}
\begin{figure}[htbp]
    \centering
    \begin{tikzpicture}[dot/.style = {circle, fill, minimum size=3pt, 
              inner sep=0pt, outer sep=0pt}  ]
        \draw plot [smooth] coordinates {(0,0) (0.4,0.6) (1,1)};
        \draw plot [smooth] coordinates {(0,0) (0.6,0.4) (1,1)};
        \draw plot [smooth] coordinates {(0,0) (0.6,-0.4) (1,-1)};
        \draw plot [smooth] coordinates {(0,0) (0.4,-0.6) (1,-1)};
        \draw plot [smooth] coordinates {(1,1) (1.7,0.85) (2.4,1)};
        \draw plot [smooth] coordinates {(1,1) (1.7,1.15) (2.4,1)};
        \draw plot [smooth] coordinates {(1,-1) (1.7,-0.85) (2.4,-1)};
        \draw plot [smooth] coordinates {(1,-1) (1.7,-1.15) (2.4,-1)};
        \draw plot [smooth] coordinates {(2.4,1) (3.1,0.85) (3.8,1)};
        \draw plot [smooth] coordinates {(2.4,1) (3.1,1.15) (3.8,1)};
        \draw plot [smooth] coordinates {(2.4,-1) (3.1,-0.85) (3.8,-1)};
        \draw plot [smooth] coordinates {(2.4,-1) (3.1,-1.15) (3.8,-1)};

        \draw plot [smooth] coordinates {(4.8,1) (5.5,0.85) (6.2,1)};
        \draw plot [smooth] coordinates {(4.8,1) (5.5,1.15) (6.2,1)};
        \draw plot [smooth] coordinates {(4.8,-1) (5.5,-0.85) (6.2,-1)};
        \draw plot [smooth] coordinates {(4.8,-1) (5.5,-1.15) (6.2,-1)};
        \draw plot [smooth] coordinates {(6.2,1) (6.9,0.85) (7.6,1)};
        \draw plot [smooth] coordinates {(6.2,1) (6.9,1.15) (7.6,1)};
        \draw plot [smooth] coordinates {(6.2,-1) (6.9,-0.85) (7.6,-1)};
        \draw plot [smooth] coordinates {(6.2,-1) (6.9,-1.15) (7.6,-1)};
        \draw plot [smooth] coordinates {(7.6,1) (8.2,0.6) (8.6,0)};
        \draw plot [smooth] coordinates {(7.6,1) (8,0.4) (8.6,0)};
        \draw plot [smooth] coordinates {(7.6,-1) (8.2,-0.6) (8.6,0)};
        \draw plot [smooth] coordinates {(7.6,-1) (8,-0.4) (8.6,0)};
        
        \draw[dotted,ultra thick] (4.15,1)--(4.45,1) (4.15,-1)--(4.45,-1);
        
        \node [dot] at (0,0) [label=left: $u$] {};
        \node [dot] at (1,1) [label=above: $x_{1,1}$] {};
        \node [dot] at (1,-1) [label=below: $x_{2,1}$] {};
        \node [dot] at (2.4,1) [label=above: $x_{1,2}$] {};
        \node [dot] at (2.4,-1) [label=below: $x_{2,2}$] {};
        \node [dot] at (3.8,1) [label={[shift={(-0.1,0)}]$x_{1,3}$}] {};
        \node [dot] at (3.8,-1) [label={[shift={(-0.1,-0.65)}]$x_{2,3}$}] {};
        \node [dot] at (4.8,1) [label={[shift={(0.1,0)}]$x_{1,a_1-2}$}] {};
        \node [dot] at (4.8,-1) [label={[shift={(0.1,-0.65)}]$x_{2,a_2-2}$}] {};
        \node [dot] at (6.2,1) [label=above: $x_{1,a_1-1}$] {};
        \node [dot] at (6.2,-1) [label=below: $x_{2,a_2-1}$] {};
        \node [dot] at (7.6,1) [label=above: $x_{1,a_1}$] {};
        \node [dot] at (7.6,-1) [label=below: $x_{2,a_2}$] {};
        \node [dot] at (8.6,0) [label=right: $v$] {};
    \end{tikzpicture}
    \caption{The schematic for the paths described in \cref{kdisjointpaths} in the case of $k = 2$. The black curves represent between vertices are paths in $G$ between those two vertices.}
    \label{kdisjointpathsimage}
\end{figure}
For a more intuitive understanding of \cref{kdisjointpaths}, it is as if there are $k$ ``meta'' paths from $u$ to $v$ where each step represents $k$ disjoint paths. \cref{kdisjointpathsimage} depicts these paths for $k = 2.$
We now consider $\FS(X,K_n)$ and show \cref{knconnectivity}, which we restate below for convenience.
\begin{customthm}{1.3}
    For $n\geq 3$, if $G = \FS(X,K_n)$ where $X$ is connected, its connectivity is equal to its minimum degree.
\end{customthm}
\begin{proof}
    We will use the fact that $\delta(G) = \lvert E(X) \rvert$ implicitly throughout the proof. We utilize the idea of atomic parts. For this argument, we will consider the atomic parts of $G.$ Recall that for any atomic part $A,$ we have $N(V(A))$ is a cutset in $\mc{C}.$ Thus, if we show that $|N(V(A))| \geq \delta(G)$ for any atomic part, this proves $G$ is $\delta(G)$-connected because the connectivity of $G$ equals $\lvert N(V(A)) \rvert\geq \delta(G)$ and is less than or equal to $\delta(G)$.    
    The general idea will be to find a structure for these atomic parts and then show that any one atomic part is adjacent to vertices of many other atomic parts. Because of \cref{adjacentatomic parts}, we then have a good lower bound on the size of the neighborhood of the atoms which we will show is larger than $\delta(G).$ 
    
    Notice that as a result of vertex transitivity, for any two vertices $\sigma$ and $\tau$ in $V(G),$ there exists an automorphism $f$ on $G$ where $f(\sigma) = \tau,$ i.e., $f(\pi) = \tau \circ\sigma^{-1}\circ \pi.$ As a result, we can group all vertices in $V(G)$ into different atomic parts $A_1, A_2, \ldots, A_w$ for some integer $w.$ These atomic parts must be disjoint by \cref{disjointatomic parts}, so they partition $V(G)$. Let $A$ and $A'$ be two of these atomic parts.
    
    Fix a vertex $\sigma$ in $V(A).$ Consider the set of transpositions $R$ given by $(i',j') \in E(X)$ such that $\sigma \circ (i',j')$ is a vertex in $V(A).$
    Let us show that for any $\tau \in V(A'),$ the vertex $\tau' = \tau \circ (i',j')$ where $(i',j') \in R$ is in $V(A')$ as well for any atomic part $A'$. 
    
    Assume for the sake of contradiction that $\tau'$ is not in $V(A').$ Then, notice that we can create a new atomic part including both $\tau$ and $\tau'$ because of the $f$ on $G$ mapping $\sigma$ to $\tau.$ As a result, $f(\sigma) = \tau$ and $f(\sigma \circ (i',j')) = \tau',$ and because there exists an atomic part containing $\sigma$ and $\sigma \circ (i',j'),$ through the automorphism, there must exist an atom containing $\tau$ and $\tau'.$ 
    But that would mean that there exist two atomic parts which are not disjoint which is a contradiction. 
    It also follows that if $\sigma \circ (k',\ell')$ is not a vertex in $A$ for some $(k',\ell') \in E(X)$ then for $\tau \in V(A')$, we have $\tau \circ (k',\ell')$ is not a vertex in $A'$ either. 
    
    Let $X'$ be the subgraph of $X$ where $V(X') = V(X)$ and $E(X') = R.$
    Notice that $A$ is isomorphic to the connected components of $\FS(X',K_n)$ since for every vertex, $\sigma \in V(\FS(X',K_n)),$ the swaps given by $(i',j') \in R$ satisfy $\sigma \circ (i',j')$ is in the same connected component of $\sigma$ similar to how the same swap of a vertex in $A$ takes the vertex to another vertex in $A.$ Therefore, each connected component is isomorphic to $A.$
    
    Let $T$ be the set of edges in $E(X)$ which have endpoints in different connected components of $X'.$ Note that $T$ and $R$ are disjoint.

    Let us first determine $\lvert V(A) \rvert.$ Because it is the size of a connected component in $\FS(X',K_n),$ we just need to determine all permutations that can be reached from a permutation, $\sigma \in \mathfrak{S}_n.$ Let the connected components of $X'$ be labeled $X_1', X_2', \ldots, X_m'.$ For each connected component, we can permute the people standing on those vertices in any way. Notice that people can't change connected components since there are no edges between components.
    Thus, the size of each component is
    $$\lvert X_1' \rvert ! \lvert X_2' \rvert ! \cdots \lvert X_m' \rvert ! = \lvert V(A) \rvert.$$

    Let us now find a lower bound on the number of atomic parts that contain a vertex adjacent to a $\sigma$ in $A.$ 
    Assume that we performed the swap $(i',j')$ to move from $\sigma$ to $\tau$ and we performed the swap $(k', \ell')$ to move from $\sigma$ to $\rho$ where $(i',j'), (k',\ell') \in T$ where $(i',j')$ and $(k',\ell')$ are different swaps.
    Let us show that $\tau$ and $\rho$ are in different atomic parts. Assume that they were in the same atomic part; then there would exist a series of swaps, $(x_1', x_2') \circ (x_3', x_4') \circ \cdots \circ (x_{c-1}', x_c')$, where $(x_{2i-1}',x_{2i}') \in R.$
    Therefore,
    $$(i',j') \circ (x_1', x_2') \circ (x_3', x_4') \circ \cdots \circ (x_{c-1}', x_c') \circ (k',\ell')$$
    is the identity permutation.

    Consider the idea of people standing on vertices in $V(X').$ If we swap $(i',j'),$ then person $\sigma(i')$ who was originally in some connected component $X_p'$ moves a different connected component $X_q',$ and the person $\sigma(j')$ who was originally in connected component $X_q'$ moves to the connected component $X_p'.$
    The swaps of the form $(x_i', x_{i+1}')$ are only swaps within connected components. In order to return the people back to their original locations, the swap $(k',\ell')$ must be returning person $\sigma(i')$ back to connected component $X_p'$ and person $\sigma(j')$ to connected component $X_q'$ but because it returns then to their exact locations that would mean that $(k',\ell') = (i',j').$ Therefore, if $(k',\ell') \neq (i',j'),$ then $\tau$ and $\rho$ must be in different atomic parts. 
    
    Thus, the number of atomic parts adjacent to a vertex in $A$ is at least $\lvert T \lvert.$ Moreover, because if a vertex in $A$ is adjacent to a vertex in $A'$ then every vertex in $A'$ is adjacent to some vertex in $A$ by \cref{adjacentatomic parts} we have that
    $$\lvert N(V(A))\rvert \geq \lvert T \lvert \cdot \lvert V(A) \rvert = \lvert T \lvert \lvert X_1' \lvert ! \lvert X_2' \lvert ! \cdots \lvert X_m' \lvert !.$$
    We will now show that $\lvert T \lvert \cdot \lvert V(A) \rvert \geq \delta(G)$ which will finish the proof.
    We will start by bounding $\delta(G).$ The maximum number of edges of $X$ within components occurs when all components are cliques.
    The number of edges between connected components of $X$ is $\lvert T \lvert.$ We know that the $\delta(G)$ is equal to the number of edges in $X.$ Let $r_1, r_2, \ldots, r_{m'}$ be a list of sizes of all components of $X'$ of size at least $2.$ 
    We have
    $$\delta(G) \leq \lvert T \rvert + \frac{r_1(r_1-1)}{2} + \frac{r_2(r_2-1)}{2} + \ldots + \frac{r_{m'}(r_{m'}-1)}{2}.$$
    Thus, because
    $$\lvert T \lvert (r_1! r_2 ! \cdots r_{m'}!-1) \geq r_1! r_2 ! \cdots r_{m'}!-1,$$
    as $\lvert T \lvert \geq 1$ since there must be multiple connected components in $X'$ in order to ensure the atoms of $G$ aren't the whole graphs themselves, it suffices to show
    $$r_1! r_2 ! \cdots r_{m'}!-1 \geq \frac{r_1(r_1-1)}{2} + \frac{r_2(r_2-1)}{2} + \ldots + \frac{r_{m'}(r_{m'}-1)}{2}.$$
    We will prove this through induction on $r_1 + \ldots + r_{m'}$ for fixed $m'$ where the base case is $r_1=\cdots=r_{m'}=2.$
    In the base case,
    $$2^{m'}-1 \geq m',$$
    which is true for all nonnegative integers $m'.$ For the inductive case assume that we have shown this for $r_1$ through $r_{m'}$ and without loss of generality we increment $r_1$ by $1.$ 
    \begin{align*}
        (r_1+1)! r_2! \cdots r_{m'}! - 1 &\geq (r_1! r_2! \cdots r_{m'}! - 1)(r_1+1),\\
        &\geq (r_1+1)\left( \frac{r_1(r_1-1)}{2} + \frac{r_2(r_2-1)}{2} + \ldots + \frac{r_{m'}(r_{m'}-1)}{2} \right),\\
        &\geq \frac{r_1+1}{r_1-1} \left( \frac{r_1(r_1-1)}{2}\right) + \frac{r_2(r_2-1)}{2} + \ldots + \frac{r_{m'}(r_{m'}-1)}{2},\\
        &\geq \frac{r_1(r_1+1)}{2} + \frac{r_2(r_2-1)}{2} + \ldots + \frac{r_{m'}(r_{m'}-1)}{2},
    \end{align*}
    finishing the proof.
\end{proof}
Let $X$ and $Y$ be graphs on $n$ vertices. We now show \cref{exchangeability} will give conditions for which $\FS(X,Y)$ is $k$-connected. We restate it below for convenience.
\begin{customprop}{1.4}
    Let $n \geq k+1.$ Assume that for any $\sigma \in V(\FS(X,Y))$ and $u',v' \in V(X)$ where $(u',v') \in E(X),$ there exist $k$ disjoint paths between $\sigma$ and $\sigma \circ (u',v').$ Then, $\FS(X,Y)$ is $k$-connected.
\end{customprop}
\begin{proof}
    Consider $\FS(X,K_n).$ The degree of each vertex is $\lvert E(X) \rvert \geq n-1$ because $X$ must be connected so that $\FS(X,K_n)$ is connected. If $\FS(X,K_n)$ is not connected, it is impossible for $\FS(X,Y)$ to be connected as $\FS(X,Y)$ is a subgraph of $\FS(X,K_n).$
    Because $n \geq k+1,$ we therefore have that $\FS(X,K_n)$ is $k$-connected by \cref{knconnectivity}.
    
    Thus between any permutations $\sigma,\rho \in V(\FS(X,K_n))$ there exist $k$ disjoint paths in $\FS(X,K_n)$. Because the adjacencies in each of these paths are given by swaps of the form $(u',v') \in E(X)$ for some intermediate permutation $\tau,$ this implies by \cref{kdisjointpaths} that $\FS(X,Y)$ is $k$-connected since the disjoint paths in $\FS(X,K_n)$ act like the metapaths of \cref{kdisjointpaths} and each of the swaps along of the metapaths can be broken down $k$-disjoint paths of $\FS(X,Y)$.
\end{proof}

\section{Connectivity of $\FS(X,\Star_n)$}\label{section:wgen2}
In this section, we prove \cref{genwilson}, which we restate below for convenience, as well as some variations and corollaries.
\begin{customthm}{1.5}
    For $n\geq 3$, if $G = \FS(X,\Star_n)$ where $X$ is connected, the connectivity of its connected components is equal to its minimum degree.
\end{customthm}
\begin{proof}
First note that $\delta(G)=\delta(X)$; we will use this fact implicitly throughout the proof. Let $k$ be the connectivity of $G.$
We begin by defining types.
\begin{definition}
    The \textit{type} of $\sigma \in V(G)$ is $\sigma^{-1}(n) \in V(X).$
\end{definition}
Types are useful because if $\sigma$ and $\rho$ in $\FS(X,\Star_n)$ are vertices of the same type, then they are essentially identical since all vertices of $\Star_n$ except the center vertex, $n,$ are essentially identical, so any permutation of those friends standing on the vertices of $X$ are essentially identical. In other words, there exists a graph automorphism $f: V(\FS(X,\Star_n)) \rightarrow V(\FS(X,\Star_n)),$ such that $f(\sigma) = \rho.$

For this argument, we will consider the atomic parts of $G.$ As mentioned at the start of \cref{knconnectivity}, it suffices to show that $k = \lvert N(V(A)) \rvert \geq \delta(G)$ for an atomic part $A.$

Let $A$ be an atomic part of $G.$ We will only be considering the atomic part $A$ from now on.

Notice that if $\lvert A \rvert \neq 1$ and there is some vertex $\sigma$ in $N(V(A))$ which is adjacent to only one vertex $\tau$ in $A,$ then $\tau$ can be removed from $A$ to get a new induced subgraph $A'.$ Therefore, $\lvert N(V(A')) \rvert \leq \lvert N(V(A)) \rvert,$ and $A'$ has less vertices than $A,$ meaning the connected components of $A'$ contain less vertices, and their neighborhoods contain at most the same number of vertices in the neighborhood of $A$, so $A$ cannot be an atomic part which is a contradiction. As a result, either $\lvert A \rvert = 1$ or all vertices in $N(V(A))$ are adjacent to at least $2$ vertices in $A.$

If $\lvert A \rvert = 1,$ we are done because $\lvert N(V(A)) \rvert = \lvert N(\sigma) \rvert \geq \delta(G)$ where $\sigma$ is the only vertex in $A.$ Therefore, we may assume all vertices in $N(V(A))$ are adjacent to at least $2$ vertices in $A.$

Let $S \subseteq V(X)$ be the set of all types included in $A.$ Assume that $\lvert S \rvert < n.$ Therefore, $N(S)$ is nonempty as $X$ is connected, so there exists some $c' \in N(S).$ Let us show that vertex $\sigma$ of type $c'$ in $N(V(A))$ cannot be adjacent to multiple vertices of $A.$ Assume for the sake of contradiction that $\rho$ of type $a'$ and $\tau$ of type $b'$ are vertices in $A$ adjacent to $\sigma.$ Notice that if $\rho$ and $\tau$ are vertices of the same type, then this is obviously true since $N(\rho)$ and $N(\tau)$ are disjoint since no vertex can neighbor two vertices of the same type due to the uniqueness of the transposition to get between two vertices of differing types. As a result, we have $\tau = \rho \circ (a',c') \circ (c',b'),$ so $\rho(c') \neq \tau(c').$

At the same time, there must exist a path between $\rho$ and $\tau$ within $A$ since $A$ is connected. Thus $\tau = \rho \circ (a',x_1') \circ(x_1',x_2') \circ \cdots \circ (x_r', b')$ for some $x_1', x_2', \ldots, x_r'.$ But notice that $x_1', x_2', \ldots, x_r'$ are all types of vertices that are contained in $A.$ Thus, there is no swap of the form $(i',j')$ where one of $i'$ or $j'$ is $c.$ Therefore, $\rho(c') = \tau(c')$ which is a contradiction.

But this contradicts that all vertices in $N(V(A))$ are adjacent to at least $2$ vertices in $A,$ so $\lvert S \rvert = n.$

Because $A$ is an atomic part meaning its neighborhood is nonempty, there exists some $\pi$ in $N(V(A)).$ Notice that there must be $\rho \in A$ where the types of $\rho$ and $\pi$ are the same as all types are contained in $S.$ Because there exists a graph isomorphism from $\rho$ to $\pi,$ there must exist an atomic part $A'$ containing $\pi.$ But then by \cref{adjacentatomic parts}, $A'$ is contained in $N(V(A)),$ so $\lvert N(V(A)) \rvert \geq n.$ Because $\delta(G) = \delta(X) < n,$ we therefore have that $\lvert N(V(A)) \rvert \geq \delta(G),$ completing the proof.
\end{proof}

Combining \cref{genwilson} with Wilson's Theorem, we have the following theorem.

\begin{theorem}\label{starkconnectivity}
    Let $X$ be a graph on $n$ vertices. Then, $\FS(X,\Star_n)$ is $k$-connected if the following conditions are satisfied:
    \begin{itemize}
        \item $X$ is biconnected,
        \item $X$ is not bipartite,
        \item $\delta(X)\geq k,$
        \item $X$ is not isomorphic to $C_n$ for $n\geq 4$,
        \item $X$ is not isomorphic to the graph $\theta_0$ on $7$ vertices shown in \cref{fig:theta0}.
    \end{itemize}
\end{theorem}

We now consider $\Star^+_n$ which will be useful in proving later results as in cases when $\FS(X,\Star_n)$ is disconnected because $X$ is bipartite, $\FS(X,\Star^+_n)$. Recall that $\Star^+_n =  \{ (i,n) : 1 \leq i \leq n-1 \} \cup \{(1,2)\}.$
We start by proving the following lemma for the result that follows.
\begin{lemma}\label{kdisjointpathsstarish}
    Let $G$ be a $k$-connected graph. Then for vertex $u$ and distinct vertices $v_i$ for $1 \leq i \leq k$ there exist disjoint paths from $u$ to $v_i$ for $1\leq i \leq k.$
\end{lemma}
\begin{proof}
    Add a new vertex $v$ to the graph which is adjacent to only $v_1,v_2,\dots,v_k$ to form a new graph, $G'.$ Notice that $G'$ is still $k$-connected. To see this, let us remove any $k-1$ vertices from $G'.$ If we remove $v,$ we must remove $k-2$ vertices from $G,$ but the graph must remain connected because $G$ is $k$-connected. If we remove $k-1$ vertices from $G,$ then $G$ must be connected, and $v$ must still be connected to the rest of $G$ because $v$ is connected to $k$ vertices in $G.$

    As a result, there exist $k$ disjoint paths from $u$ to $v,$ but this also generates disjoint paths from $u$ to $v_i$ for $1 \leq i \leq k.$
\end{proof}
We can now determine when $\FS(\Star^+_n,Y)$ is $k$-connected.
\begin{theorem}\label{star+kconnectivity}
    Let $X$ be a graph on $n \geq 3$ vertices. Then, $\FS(X,\Star^+_n)$ is $k$-connected if the following conditions are satisfied:
    \begin{itemize}
        \item $X$ is biconnected,
        \item $\delta(X) \geq k,$
        \item $X$ is not isomorphic to $C_n$ for $n\geq 4$,
        \item $X$ is not isomorphic to the graph $\theta_0.$
    \end{itemize}
\end{theorem}
\begin{proof}
    Notice that if $X$ is not bipartite, then $\FS(X, \Star^+_n)$ is $k$-connected because $\FS(X, \Star_n)$ is $k$-connected. Therefore, let $X$ be bipartite with a bipartition of $A$ and $B$. As mentioned in the Remark 2.8 by Defant and Kravitz \cite{defant2021friends}, the connected components of $\FS(X, \Star_n)$ are given by
    $$\{\sigma \in \mathfrak{S}_n : p(\sigma) \equiv 0 \pmod{2}\} \qquad \qquad \{\sigma \in \mathfrak{S}_n : p(\sigma) \equiv 1 \pmod{2}\},$$
    where
    $$p(\sigma) = \lvert \sigma([n-1])\cap A \rvert + \frac{\sgn(\sigma)+1}{2}.$$
    Because $X$ is biconnected, we know that the two connected components are both $k$-connected by \cref{genwilson} since $\delta(X) \geq k$. 
    We now consider the additional edges created in $\FS(X,\Star^+_n)$ as a result of the edge $(1,2) \in E(\Star^+_n).$ Let $\sigma \in \mathfrak{S}_n$ be some permutation. Consider $\sigma' = \sigma \circ (1,2).$
    Notice that
    $$p(\sigma) \equiv \lvert \sigma([n-1])\cap A \rvert + \frac{\sgn(\sigma)+1}{2} \equiv \lvert \sigma'([n-1])\cap A \rvert + \left(\frac{\sgn(\sigma')+1}{2} + 1 \right) \equiv p(\sigma')+1 \pmod{2}.$$
    Therefore, swapping along $(1,2)$ creates an edge from the first connected component to the second connected component. Because $\delta(Y)\geq k$, we have $\lvert E(Y) \rvert \geq k,$ so there are at least $k$ edges formed by swapping along $(1,2).$
    Let these edges be between $\rho_i \in \mathfrak{S}_n$ in the first connected component and $\rho_i' \in \mathfrak{S}_n$ in the second connected component for $1\leq i \leq k.$

    Consider any two vertices, $\tau_1, \tau_2 \in \mathfrak{S}_n.$ If $\tau_1$ and $\tau_2$ are in the same connected component, then there already exist $k$ disjoint paths between the two vertices since each component is $k$-connected.
    If $\tau_1$ and $\tau_2$ are in different components, say the first and second component, respectively, there exist disjoint paths from $\tau_1$ to $\rho_i$ for $1 \leq i \leq k$ by \cref{kdisjointpathsstarish}.
    From there, the edge swap along $(1,2)$ converts $\rho_i$ into $\rho_i'$ for $1\leq i \leq k.$ Finally, there exist disjoint paths from $\rho_i'$ for $1\leq i \leq k$ to $\tau_2$ by \cref{kdisjointpathsstarish}, finishing the proof that $\FS(X,\Star^+_n)$ is $k$-connected.  
\end{proof}

Note that the last two conditions in both \cref{starkconnectivity,star+kconnectivity} are irrelevant for $k\geq3$.

\section{Bounds on Minimum Degree}\label{section:bangachev}
\subsection{\cref{kban1}}
This subsection is devoted to proving \cref{kban1}, which we restate below for convenience.
\begin{customthm}{1.6}
    Let $X$ and $Y$ be two graphs on $n$ vertices satisfying
    \begin{itemize}
        \item $\delta(X), \delta(Y)>n/2,$
        \item $2\min\{\delta(X), \delta(Y)\}+3\max\{\delta(X),\delta(Y)\} \geq 3n+2k-4,$
    \end{itemize}
    where $k \geq 2.$ Then for sufficiently large $n$, namely $n \geq 5(1+(k-1)(k+6)+2k-3),$ we have $\FS(X,Y)$ is $k$-connected.
\end{customthm}
First note that $n\geq 5(1+(k-1)(k+6)+2k-3) \geq 2k+6$ for $k\geq 2.$ This will be used implicitly throughout this section.

Define $Y$ to be \textit{$k$-Wilsonian} to mean that $\FS(\Star_n,Y)$ is $k$-connected, and define $Y$ to be \textit{almost $k$-Wilsonian} if $\FS(\Star^+_n,Y)$ is $k$-connected. We then have the following lemma. Essentially the same proof as the one in Bangachev's paper \cite{bangachev2022asymmetric} holds in vertex connectivity where we use \cref{starkconnectivity} and \cref{star+kconnectivity} to deduce claims involving vertex connectivity, so the proof of the lemma will be omitted.
\begin{lemma}[\protect{\cite[Lemma 4.1]{bangachev2022asymmetric}}]
    Let $G$ be a graph on $m$ vertices with minimum degree $\delta(G)$. Let $Q$ be any subset of $V(G)$ such that we have $\lvert Q \rvert \geq k+5$ and we have $2\lvert Q\rvert + 3\delta(G) \geq 3m +2k- 2.$ Then $G\lvert_Q$ has at most two connected components.
    Furthermore:
    \begin{enumerate}
        \item If $G\lvert_Q$ has exactly two connected components, $F_1$ and $F_2$, then both $F_1$ and $F_2$ are $k$-Wilsonian, and the following inequalities are satisfied:
        \begin{align*}
            &\delta(G)+1+\lvert Q\rvert -m \leq \lvert V(F_1) \rvert \leq m-\delta(G)-1, \quad \delta(F_1) \geq \delta(G) + \lvert Q \rvert - m,\\
            &\delta(G)+1+\lvert Q\rvert -m \leq \lvert V(F_2) \rvert \leq m-\delta(G)-1, \quad \delta(F_2) \geq \delta(G) + \lvert Q \rvert - m.
        \end{align*}
        \item If $G\lvert_Q$ has a single connected component $F$, then one of the following holds:
        \begin{itemize}
            \item $F$ is almost $k$-Wilsonian
            \item There exists a cut vertex, $v,$ such that $F\lvert_{V(F)\backslash \{v\}}$ has exactly two connected components, $F_1'$ and $F_2'.$ Furthermore both $F_1'$ and $F_2'$ are $k$-Wilsonian, and the following inequalities are satisfied:
            \begin{align*}
                &\delta(G)+\lvert Q \rvert - m \leq \lvert V(F_1') \rvert \leq m-\delta(G)-1, \quad \delta(F_1') \geq \delta(G) + \lvert Q \rvert - m - 1,\\
                &\delta(G)+\lvert Q \rvert - m \leq \lvert V(F_2') \rvert \leq m-\delta(G)-1, \quad \delta(F_2') \geq \delta(G) + \lvert Q \rvert - m - 1. 
            \end{align*}
            Finally if $\lvert N(v) \cap V(F_1') \rvert \geq k,$ then the graph $F\lvert_{V(F_1')\cup \{v\}}$ is $k$-Wilsonian, and likewise for $F_2'.$
        \end{itemize}
    \end{enumerate}
\end{lemma}

Without loss of generality, let $\delta(X) \leq \delta(Y)$ so that the second condition on the minimum degree is $2\delta(X)+3\delta(Y)\geq 3n+2k-4.$

Let us start by fixing a bijection $\sigma.$ Let $u,v \in V(Y)$ such that $(u',v') \in E(X)$ where $u' = \sigma^{-1}(u)$ and $v' = \sigma^{-1}(v).$ If we show that there exist $k$ disjoint paths between $\sigma$ and $\sigma\circ (u',v'),$ then because $n \geq k+1$ and $X$ is connected as $\delta(X) > n/2$, we must have by \cref{exchangeability} that $\FS(X,Y)$ is $k$-connected.
Assume for the sake of contradiction that there do not exist $k$ disjoint paths between $\sigma$ and $\sigma \circ (u',v').$

For an edge, $e = (i',j') \in E(X),$ define $\sigma(e) = (i,j)$ where $i = \sigma(i')$ and $j = \sigma(j').$ Let us first consider the case when $u$ and $v$ are adjacent. Consider the edges $(i,j) \in E(Y)$ where $(i',j') \in E(X)$ as well. Note that $n\geq 2(k+12)$ for $k\geq 2.$ We have
\begin{align*}
    \lvert E(Y) \cap \sigma(E(X)) \lvert &\geq \frac{n\delta(X)}{2} + \frac{n\delta(Y)}{2} - \frac{n(n-1)}{2},\\
    &= \frac{n}{2}\left( \delta(X)+\delta(Y)-n+1\right),\\
    &=  \frac{n}{2}\left(\frac{2\delta(X)+3\delta(Y)}{3}+\frac{\delta(X)}{3}-n+1\right)\\
    &> \frac{n}{2}\left(\frac{3n+2k-4}{3}+\frac{k+12}{3}-n+1\right) \geq 3n \geq 2n-3+(k-1).
\end{align*}
As a result there must exist $k-1$ such edges where $i \not \in \{u,v\}$ and $j \not \in \{u,v \}$, as there are at most $2n-3$ edges incident to either $u$ or $v$.
So, our paths can be to first swap along that edge, then swap $(u,v)$ then swap along the edge again to create $k-1$ disjoint paths. Alternatively, we directly swap along $(u,v)$ creating one more disjoint path, completing this case of the proof of \cref{kban1}.

We now assume that $(u,v)\not\in E(Y).$ Let $A' = N(u') \cap N(v'),$ $B = N(u) \cap N(v),$ $A = \sigma(A'),$ and $B' = \sigma^{-1}(B).$ The following lemmas are slight modifications of lemmas in Bangachev \cite{bangachev2022asymmetric} and so their proofs are not included. The lemmas assume that there do not exist $k$-disjoint paths between $\sigma$ and $\sigma \circ (u',v'),$ and in the proofs of Bangachev \cite{bangachev2022asymmetric}, exchangeability is replaced with $k$-exchangeability which refers to the existence of $k$-disjoint paths between $\sigma$ and $\sigma\circ (u',v')$ for all $\sigma.$
\begin{lemma}[\protect{\cite[Lemma 4.2]{bangachev2022asymmetric}}]
    The following inequalities hold:
    \begin{align*}
        \lvert B \cap \sigma(N(u')) \lvert &\leq k-1,\\
        \lvert B \cap \sigma(N(v')) \lvert &\leq k-1.
    \end{align*}
\end{lemma}
\begin{lemma}[\protect{\cite[Lemma 4.3]{bangachev2022asymmetric}}]
    The following inequalities hold:
    \begin{align*}
        \lvert B \lvert &\geq 2 \delta(Y) + 2-n,\\
        \lvert A' \lvert &\geq 2\delta(X) + 2\delta(Y) - 2n-2k+4.
    \end{align*}
\end{lemma}

Notice that because $\delta(Y) \geq \delta(X),$ we have that $\delta(Y) \geq 3n/5.$ Therefore, we have
\[ \lvert B \rvert \geq 2\delta(Y) + 2 - n \geq n/5 + 2. \]
Let $W = B \setminus (\sigma(N(u')) \cup \sigma(N(v'))).$ We also know that $\lvert B \cap \sigma(N(u')) \rvert \leq k-1$ and $\lvert B \cap \sigma(N(v')) \rvert \leq k-1,$ so 
$$\lvert W \rvert \geq n/5+4-2k \geq 1+(k-1)(k+6).$$ 

Let $w \in W.$ Denote $w' = \sigma^{-1}(w).$
Consider the intersection between $A$ and $N(w).$ Note that $w \not \in A$ because $w \not \in \sigma(N(u'))$.
We have
$$\lvert N(w) \cap A \lvert =  \lvert N[w] \cap A \lvert \geq \lvert N[w] \rvert + \lvert A \rvert -n \geq (\delta(Y) + 1) + (2\delta(X) + 2\delta(Y) - 2n-2k+4) - n \geq 1.$$
Thus, $w$ has a neighbor $x(w) \in A.$ Let $x'(w) = \sigma^{-1}(x(w)).$
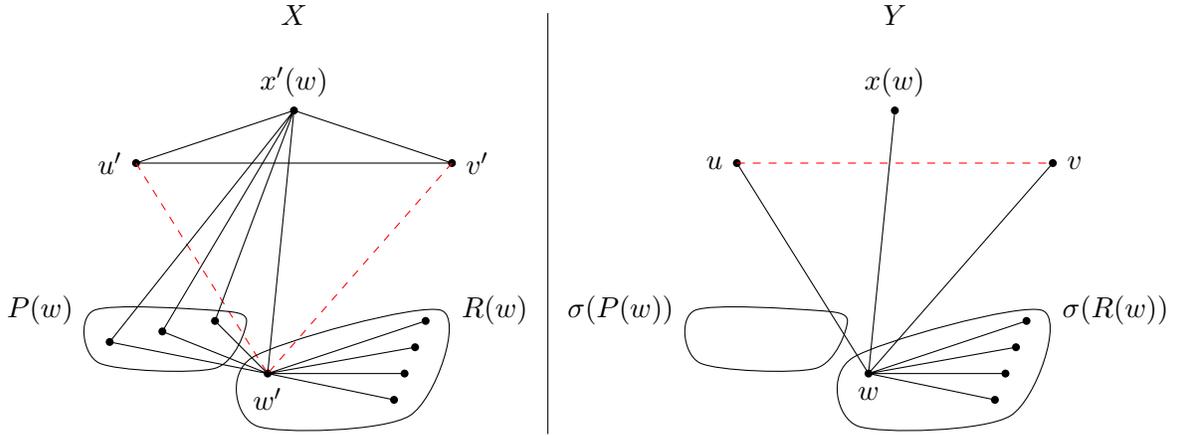
\begin{figure}[htbp]
    \begin{tikzpicture}[scale=0.7,
        dot/.style = {circle, fill, minimum size=3pt, 
              inner sep=0pt, outer sep=0pt}    
        ]
        \coordinate (w') at (0,0);
        \coordinate (x') at (0.5,5);
        \coordinate (u') at (-2.5,4);
        \coordinate (v') at (3.5,4);
        \node [dot] at (w') [label=below: $w'$] {};
        \node [dot] at (x') [label=above: $x'(w)$] {};
        \node [dot] at (u') [label=left: $u'$] {};
        \node [dot] at (v') [label=right: $v'$] {};
        \draw (x') node[above = 1 cm] {$X$};
        \draw (u')--(x')--(v')--(u') (x')--(w');
        \draw[dashed,red] (u')--(w') (v')--(w');
    \end{tikzpicture}
    \begin{tikzpicture}[scale=0.7]
        \draw (0,8)--(0,0);
    \end{tikzpicture}
    \begin{tikzpicture}[scale=0.7,
        dot/.style = {circle, fill, minimum size=3pt, 
              inner sep=0pt, outer sep=0pt}    
        ]
        \coordinate (w) at (0,0);
        \coordinate (x) at (0.5,5);
        \coordinate (u) at (-2.5,4);
        \coordinate (v) at (3.5,4);
        \node [dot] at (w) [label=below: $w$] {};
        \node [dot] at (x) [label=above: $x(w)$] {};
        \node [dot] at (u) [label=left: $u$] {};
        \node [dot] at (v) [label=right: $v$] {};
        \draw (x) node[above = 1 cm] {$Y$};
        \draw (w)--(u) (w)--(x) (w)--(v);
        \draw[dashed,red] (u)--(v);
    \end{tikzpicture}
    \caption{The vertices $w', u', v',$ and $x'(w)$ as well as the corresponding vertices in $Y.$ Black edges represent edges in the graph and red edges represent edges which do not exist in the graph.}
    \label{fig:XYgraph1}
\end{figure}
\begin{lemma}
There are at most $k-1$ elements $w \in W$ for which $x'(w)$ is adjacent to $w' = \sigma^{-1}(w).$
\end{lemma}
\begin{proof}
    Assume for the sake of contradiction that there are at least $k$ elements $w \in W$ which satisfy $x'(w)$ is adjacent to $w' = \sigma^{-1}(w)$ as shown in \cref{fig:XYgraph1}. Notice that for each of those elements we can construct a path,
    $$(w',x'(w)), (x'(w),u'), (u',v'), (v',x'(w)), (x'(w),w').$$
    Notice that for two distinct such $w, w^\star \in W,$ these paths are disjoint since after the first swap until the last swap, the first path has $w$ in a different position than $\sigma^{-1}(w)$ and $w^\star$ in the same position as the start while in the second path the opposite is true.
    Thus, these two paths must be disjoint.
    Thus, with these $k$ elements $w\in W$ we can construct $k$ disjoint paths between $\sigma$ and $\sigma\circ(u', v')$, contradicting our assumption that $k$ such paths do not exist.
\begin{figure}[htbp]
    \begin{tikzpicture}[scale=0.7,
        dot/.style = {circle, fill, minimum size=3pt, 
              inner sep=0pt, outer sep=0pt}    
        ]
        \coordinate (w') at (0,0);
        \coordinate (x') at (0.5,5);
        \coordinate (u') at (-2.5,4);
        \coordinate (v') at (3.5,4);
        \coordinate (P1) at (-1, 1);
        \coordinate (P2) at (-2, 0.8);
        \coordinate (P3) at (-3, 0.6);
        \coordinate (R1) at (3,1);
        \coordinate (R2) at (2.8,0.5);
        \coordinate (R3) at (2.6,0);
        \coordinate (h') at (2.4,-0.5);
        \node [dot] at (w') [label=below: $w'$] {};
        \node [dot] at (x') [label=above: $x'(w)$] {};
        \node [dot] at (u') [label=left: $u'$] {};
        \node [dot] at (v') [label=right: $v'$] {};
        \node [dot] at (P1) {};
        \node [dot] at (P2) {};
        \node [dot] at (P3) {};
        \node [dot] at (R1) {};
        \node [dot] at (R2) {};
        \node [dot] at (R3) {};
        \node [dot] at (h') [label=right: $h'(w)$] {};
        \draw (x') node[above = 1 cm] {$X$};
        \draw (h')--(x');
        \draw (u')--(x')--(v')--(u');
        \draw[dashed,red] (x')--(w') (u')--(w') (v')--(w');
        \draw (w')--(P1)--(x') (w')--(P2)--(x') (w')--(P3)--(x');
        \draw (w')--(R1) (w')--(R2) (w')--(R3); 
        \draw (w')--(h');
        \draw plot [smooth cycle] coordinates {(-1,0.1) (-0.4,1) (-1,1.2) (-3.3,1.2) (-3.2,0.2)};
        \draw plot [smooth cycle] coordinates {(-0.3,0.3) (-0.2,-1) (4,-1) (3.3,1.2)};
        \draw (-3.3,1.2) node [label=left: $P(w)$] {};
        \draw (3.3,1.2) node [label=right: $R(w)$] {};
    \end{tikzpicture}
    \begin{tikzpicture}[scale=0.7]
        \draw (0,8)--(0,0);
    \end{tikzpicture}
    \begin{tikzpicture}[scale=0.7,
        dot/.style = {circle, fill, minimum size=3pt, 
              inner sep=0pt, outer sep=0pt}    
        ]
        \coordinate (w) at (0,0);
        \coordinate (x) at (0.5,5);
        \coordinate (u) at (-2.5,4);
        \coordinate (v) at (3.5,4);
        \coordinate (P1) at (-1, 1);
        \coordinate (P2) at (-2, 0.8);
        \coordinate (P3) at (-3, 0.6);
        \coordinate (R1) at (3,1);
        \coordinate (R2) at (2.8,0.5);
        \coordinate (R3) at (2.6,0);
        \coordinate (h) at (2.4,-0.5);
        \node [dot] at (w) [label=below: $w$] {};
        \node [dot] at (x) [label=above: $x(w)$] {};
        \node [dot] at (u) [label=left: $u$] {};
        \node [dot] at (v) [label=right: $v$] {};
        \node [dot] at (R1) {};
        \node [dot] at (R2) {};
        \node [dot] at (R3) {};
        \node [dot] at (h) [label=right: $h(w)$] {};
        \draw (x) node[above = 1 cm] {$Y$};
        \draw (w)--(u) (w)--(x) (w)--(v);
        \draw[dashed,red] (u)--(v);
        \draw (w)--(R1) (w)--(R2) (w)--(R3);
        \draw (w)--(h);
        \draw plot [smooth cycle] coordinates {(-1,0.1) (-0.4,1) (-1,1.2) (-3.3,1.2) (-3.2,0.2)};
        \draw plot [smooth cycle] coordinates {(-0.3,0.3) (-0.2,-1) (4,-1) (3.3,1.2)};
        \draw (-3.3,1.2) node [label=left: $\sigma(P(w))$] {};
        \draw (3.3,1.2) node [label=right: $\sigma(R(w))$] {};
    \end{tikzpicture}
    \caption{The vertices $w', u', v', h'(w)$ and $x'(w)$ and the sets of vertices $P(w)$ and $R(w)$ as well as the corresponding vertices and sets in $Y.$ Black edges represent edges in the graph and red edges represent edges which do not exist in the graph.}
    \label{fig:XYgraph2}
\end{figure}
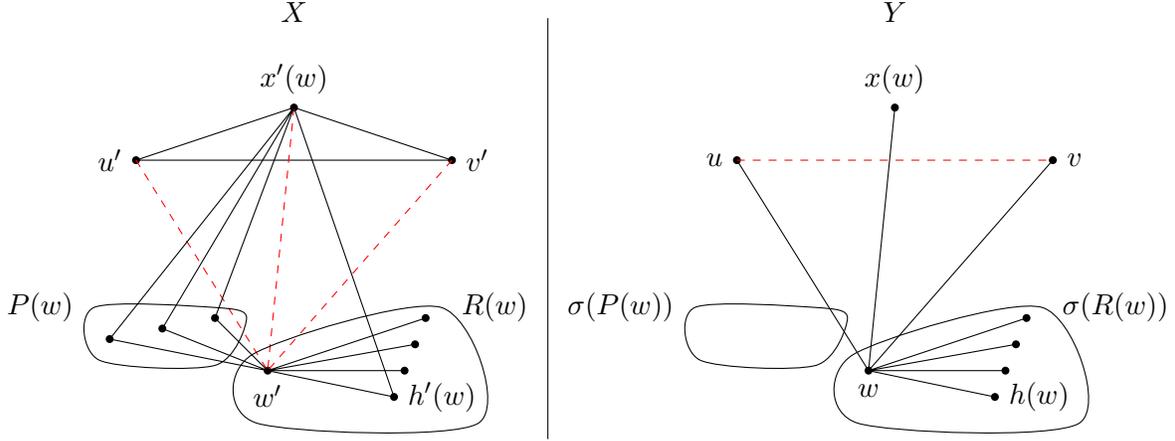
\end{proof}
Define $W_1$ to be the subset of $W$ consisting of the $w\in W$ such that $x'(w)$ is not adjacent to $w' = \sigma^{-1}(w).$ We will now consider the following sets
\begin{align*}
    P(w) &= N(x'(w)) \cap N(w'),\\
    R(w) &= N[w'] \cap \sigma^{-1}(N[w]).
\end{align*}
Note that $u', v' \not \in P(w)$ and $u', v' \not \in R(w)$ because $u', v' \not \in N(w').$
\begin{lemma}
    There exist at most $2k-2$ elements $w \in W_1$ where $R(w)$ and $N[x'(w)]$ are not disjoint.
\end{lemma}
\begin{proof}
    If $w$ satisfies that $R(w)$ and $N[x'(w)]$ share a vertex $h(w)$ as shown in \cref{fig:XYgraph2} then the $(X,Y)$-friendly swaps,
    $$(w',h'(w)), (h'(w),x'(w)),(x'(w),u'),(u',v'),(v',x'(w)),(x'(w),h'(w)),(h'(w),w'),$$
    swaps $u$ and $v$ while maintaining the positions of the other vertices. Consider $w,w^\star \in W_1.$ We refer to the path formed by the swaps above as the first path, and the path formed by the swaps where $w$ is replaced by $w^\star$ as the second path. Notice that different $w$ and $w^\star$ must have disjoint corresponding paths if $w \neq h(w^\star)$ or $w^\star \neq h(w).$ If $w \neq h(w^\star)$ and $w^\star \neq h(w),$ this is true because in the first path the position of $w$ changes while the position of $w^\star$ stays the same, and in the second path, the position of $w^\star$ changes while the position of $w$ stays the same. Therefore, in intermediate permutations of the first path, $w'$ corresponds to a person who is not $w$ and $w^{\star\prime}$ corresponds to $w^\star$ while in the second path $w'$ corresponds to $w$ and $w^{\star\prime}$ corresponds to a person who is not $w^\star$ meaning none of the permutations are shared. If $w \neq h(w^\star)$ and $w^\star = h(w),$ then after the first swap of the first path, the position of $w^\star$ is $\sigma^{-1}(w)$ and in the second path, $w^\star$ never reaches this position as $w \neq h(w^\star)$. Thus, the paths are disjoint. By similar logic if $w^\star \neq h(w^\star)$ and $w^\star = h(w),$ the paths are disjoint. It follows that there exist at most $2k-2$ elements $w \in W_1$ where $R(w)$ and $N[x'(w)]$ are not disjoint as any more elements would create $k$ disjoint paths since for each $w$ in $W_1,$ there is at most one other $w^\star$, i.e., $w^\star = h(w),$ which can have a path which shares a vertex with the path of $w$.
\end{proof}
Define $W_2$ to be the subset of $W_1$ consisting of the $w \in W_1$ such that $R(w)$ and $N[x'(w)]$ are disjoint.
\begin{lemma}
    There are at most $(k-1)(k+3)$ elements $w \in W_2$ such that some vertex in $\sigma(P(w))$ has more than one neighbor in $\sigma(R(w)).$
\end{lemma}
\begin{proof}
    Assume that $w \in W_2$ satisfies that $h(w) \in \sigma(P(w))$ has more than one neighbor in $\sigma(R(w)).$ Let $q_1(w)$ and $q_2(w)$ be these neighbors. Let $h'(w) = \sigma^{-1}(h(w)),$ $q_1'(w) = \sigma^{-1}(q_1(w))$ and $q_2'(w) = \sigma^{-1}(q_2(w)).$ We illustrate these vertices in \cref{fig:XYgraph3}.

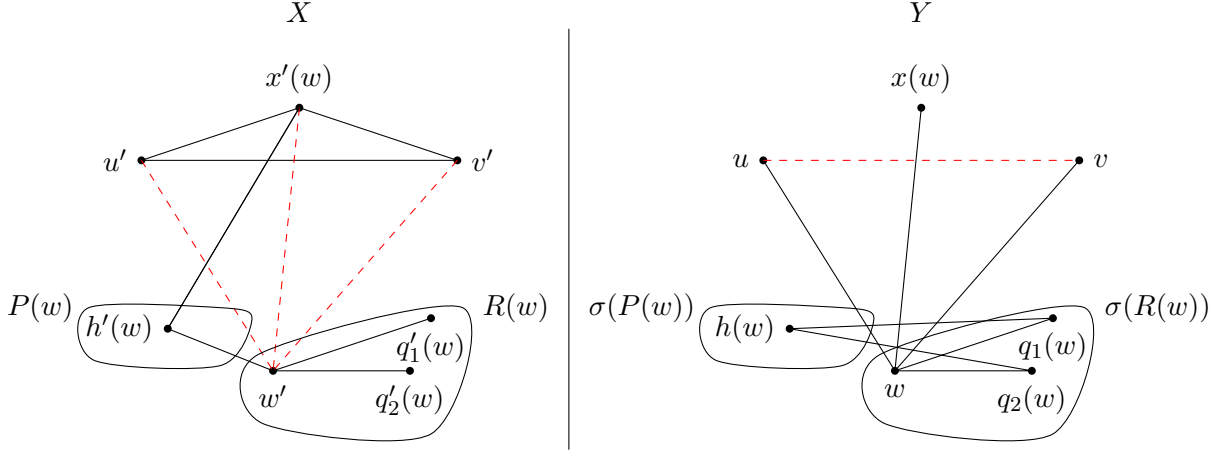
\begin{figure}[htbp]
    \begin{tikzpicture}[scale=0.7,
        dot/.style = {circle, fill, minimum size=3pt, 
              inner sep=0pt, outer sep=0pt}    
        ]
        \coordinate (w') at (0,0);
        \coordinate (x') at (0.5,5);
        \coordinate (u') at (-2.5,4);
        \coordinate (v') at (3.5,4);
        \coordinate (P1) at (-1, 1);
        \coordinate (h') at (-2, 0.8);
        \coordinate (P3) at (-3, 0.6);
        \coordinate (q1') at (3,1);
        \coordinate (R2) at (2.8,0.5);
        \coordinate (q2') at (2.6,0);
        \coordinate (R4) at (2.4,-0.5);
        \node [dot] at (w') [label=below: $w'$] {};
        \node [dot] at (x') [label=above: $x'(w)$] {};
        \node [dot] at (u') [label=left: $u'$] {};
        \node [dot] at (v') [label=right: $v'$] {};
        \node [dot] at (h') [label=left: $h'(w)$] {};
        \node [dot] at (q1') [label=below: $q_1'(w)$] {};
        \node [dot] at (q2') [label=below: $q_2'(w)$] {};
        \draw (x') node[above = 1 cm] {$X$};
        \draw (h')--(x');
        \draw (u')--(x')--(v')--(u');
        \draw[dashed,red] (x')--(w') (u')--(w') (v')--(w');
        \draw (w')--(h')--(x') (w')--(q1') (w')--(q2');
        \draw plot [smooth cycle] coordinates {(-1,0.1) (-0.4,1) (-1,1.2) (-3.5,1.2) (-3.4,0.2)};
        \draw plot [smooth cycle] coordinates {(-0.3,0.3) (-0.2,-1) (3,-1.2) (3.6,1.2)};
        \draw (-3.4,1.2) node [label=left: $P(w)$] {};
        \draw (3.6,1.2) node [label=right: $R(w)$] {};
    \end{tikzpicture}
    \begin{tikzpicture}[scale=0.7]
        \draw (0,8)--(0,0);
    \end{tikzpicture}
    \begin{tikzpicture}[scale=0.7,
        dot/.style = {circle, fill, minimum size=3pt, 
              inner sep=0pt, outer sep=0pt}    
        ]
        \coordinate (w) at (0,0);
        \coordinate (x) at (0.5,5);
        \coordinate (u) at (-2.5,4);
        \coordinate (v) at (3.5,4);
        \coordinate (P1) at (-1, 1);
        \coordinate (h) at (-2, 0.8);
        \coordinate (P3) at (-3, 0.6);
        \coordinate (q1) at (3,1);
        \coordinate (R2) at (2.8,0.5);
        \coordinate (q2) at (2.6,0);
        \coordinate (R4) at (2.4,-0.5);
        \node [dot] at (w) [label=below: $w$] {};
        \node [dot] at (x) [label=above: $x(w)$] {};
        \node [dot] at (u) [label=left: $u$] {};
        \node [dot] at (v) [label=right: $v$] {};
        \node [dot] at (h) [label=left: $h(w)$] {};
        \node [dot] at (q1) [label=below: $q_1(w)$]{};
        \node [dot] at (q2)[label=below: $q_2(w)$] {};
        \draw (x) node[above = 1 cm] {$Y$};
        \draw (w)--(u) (w)--(x) (w)--(v);
        \draw[dashed,red] (u)--(v);
        \draw (h)--(q1) (h)--(q2) (w)--(q1) (w)--(q2);
        \draw plot [smooth cycle] coordinates {(-1,0.1) (-0.4,1) (-1,1.2) (-3.5,1.2) (-3.4,0.2)};
        \draw plot [smooth cycle] coordinates {(-0.3,0.3) (-0.2,-1) (3,-1.2) (3.6,1.2)};
        \draw (-3.4,1.2) node [label=left: $\sigma(P(w))$] {};
        \draw (3.6,1.2) node [label=right: $\sigma(R(w))$] {};
    \end{tikzpicture}
    \caption{The vertices $w', u', v', h'(w), q_1'(w), q_2'(w)$ and $x'(w)$ and the sets of vertices $P(w)$ and $R(w)$ as well as the corresponding vertices and sets in $Y.$ Black edges represent edges in the graph and red edges represent edges which do not exist in the graph.}
    \label{fig:XYgraph3}
\end{figure}

    Notice that $w \not \in \{q_1(w), q_2(w)\}$ since that would imply $h(w) \in R(w)$ which is impossible by the assumption that $R(w)$ and $N[x'(w)]$ are disjoint. We therefore have the following sequence of $(X,Y)$-swaps which swaps the position of $u$ and $v.$ The path below will be referred to as the path of $w.$
    
    \begin{gather*}
    (w',q_1'(w)),(h'(w),w'),(w',q_2'(w)),(w',q_1'(w)),(w',h'(w)),(h'(w),x'(w)),(x'(w),u'),(u',v'),\\(v',x'(w)),(h'(w),x'(w)),(w',h'(w)),(w',q_1'(w)),(w',q_2'(w)),(h'(w),w'),(w',q_1'(w)).
    \end{gather*}
    Consider $w^\star \in W_2$ such that $w^\star \not \in \{ q_1(w), q_2(w), h(w) \}.$ Let us now determine all $w^\star$ whose corresponding path might share a vertex with the above path. We refer to the path created by replacing $w$ with $w^\star$ in the swaps above as the path of $w^\star.$ Consider the positions person $w^\star$ reaches. Notice that in the path of $w,$ the position of $w^\star$ remains fixed, so in the path of $w^\star,$ the only intermediate permutations which could be shared with the path of $w$ are $\tau$ where $\tau^{-1}(w^\star) = \sigma^{-1}(w^\star).$ 

    There are exactly $2$ such permutations in the entire path, $\rho_1^{(w^\star)}$ and $\rho_2^{(w^\star)},$ where $\rho_1^{(w^\star)}(h'(w^\star)) = q_1(w^\star),$ $\rho_1^{(w^\star)}(q_1'(w^\star)) = q_2(w),$ $\rho_1^{(w^\star)}(q_2'(w^\star)) = h(w^\star),$ $\rho_1^{(w^\star)}(x'(w^\star)) = x(w^\star),$ $\rho_1^{(w^\star)}(u') = u,$ and $\rho_1^{(w^\star)}(v') = v$ and $\rho_2^{(w^\star)}$ is identical except $\rho_2^{(w^\star)}(u') = v,$ and $\rho_2^{(w)}(v') = u.$ We will only consider $\rho_1^{(w^\star)}$ because if $\rho_1^{(w^\star)}$ is shared in the path of $w,$ then $\rho_2^{(w^\star)}$ must be shared as well and vice versa.

    Notice that $\rho_1^{(w^\star)}$ just cycles $h(w^\star), q_1(w^\star),$ and $q_2(w^\star).$ Thus, we just need to find all intermediate permutations, $\tau,$ in the path of $w$ where the permutation cycles the positions of three people. This occurs only at $\rho_1^{(w)}$ and $\pi$ where $\pi(h'(w)) = q_1(w),$ $\pi(w') = h(w),$ and $\pi(q_1'(w)) = w.$ These permutations are shown in \cref{fig:rhopi}.

\begin{figure}[htbp]
    \begin{tikzpicture}[scale=0.7,
        dot/.style = {circle, fill, minimum size=3pt, 
              inner sep=0pt, outer sep=0pt}    
        ]
        \coordinate (w') at (0,0);
        \coordinate (x') at (0.5,5);
        \coordinate (u') at (-2.5,4);
        \coordinate (v') at (3.5,4);
        \coordinate (P1) at (-1, 1);
        \coordinate (h') at (-2, 0.8);
        \coordinate (P3) at (-3, 0.6);
        \coordinate (q1') at (3,1);
        \coordinate (R2) at (2.8,0.5);
        \coordinate (q2') at (2.6,0);
        \coordinate (R4) at (2.4,-0.5);
        \node [dot] at (w') [label=below: $w$] {};
        \node [dot] at (x') [label=above: $x(w)$] {};
        \node [dot] at (u') [label=left: $u$] {};
        \node [dot] at (v') [label=right: $v$] {};
        \node [dot] at (h') [label=left: $q_1(w)$] {};
        \node [dot] at (q1') [label=below: $q_2(w)$] {};
        \node [dot] at (q2') [label=below: $h(w)$] {};
        \draw (x') node[above = 1 cm] {$\rho_1^{(w)}$};
        \draw (h')--(x');
        \draw (u')--(x')--(v')--(u');
        \draw (w')--(h')--(x') (w')--(q1') (w')--(q2');
        \draw plot [smooth cycle] coordinates {(-1,0.1) (-0.4,1) (-1,1.2) (-3.5,1.2) (-3.4,0.2)};
        \draw plot [smooth cycle] coordinates {(-0.3,0.3) (-0.5,-1) (3,-1.2) (3.6,1.2)};
        \draw (-3.4,1.2) node [label=left: $P(w)$] {};
        \draw (3.6,1.2) node [label=right: $R(w)$] {};
    \end{tikzpicture}
    \begin{tikzpicture}[scale=0.7]
        \draw (0,8)--(0,0);
    \end{tikzpicture}
    \begin{tikzpicture}[scale=0.7,
        dot/.style = {circle, fill, minimum size=3pt, 
              inner sep=0pt, outer sep=0pt}    
        ]
        \coordinate (w') at (0,0);
        \coordinate (x') at (0.5,5);
        \coordinate (u') at (-2.5,4);
        \coordinate (v') at (3.5,4);
        \coordinate (P1) at (-1, 1);
        \coordinate (h') at (-2, 0.8);
        \coordinate (P3) at (-3, 0.6);
        \coordinate (q1') at (3,1);
        \coordinate (R2) at (2.8,0.5);
        \coordinate (q2') at (2.6,0);
        \coordinate (R4) at (2.4,-0.5);
        \node [dot] at (w') [label=below: $h(w)$] {};
        \node [dot] at (x') [label=above: $x(w)$] {};
        \node [dot] at (u') [label=left: $u$] {};
        \node [dot] at (v') [label=right: $v$] {};
        \node [dot] at (h') [label=left: $q_1(w)$] {};
        \node [dot] at (q1') [label=below: $w$] {};
        \node [dot] at (q2') [label=below: $q_2(w)$] {};
        \draw (x') node[above = 1 cm] {$\pi$};
        \draw (h')--(x');
        \draw (u')--(x')--(v')--(u');
        \draw (w')--(h')--(x') (w')--(q1') (w')--(q2');
        \draw plot [smooth cycle] coordinates {(-1,0.1) (-0.4,1) (-1,1.2) (-3.5,1.2) (-3.4,0.2)};
        \draw plot [smooth cycle] coordinates {(-0.3,0.3) (-0.5,-1) (3,-1.2) (3.6,1.2)};
        \draw (-3.4,1.2) node [label=left: $P(w)$] {};
        \draw (3.6,1.2) node [label=right: $R(w)$] {};
    \end{tikzpicture}
    \caption{Depictions of $\rho_1^{(w)}$ and $\pi$. In each depiction, the graph of the underlying $X$ graph is shown. The labels of the vertices are the people standing on the vertices of the graph. The vertices of the $X$ graph are placed in analogous positions to \cref{fig:XYgraph3}.}
    \label{fig:rhopi}
\end{figure}
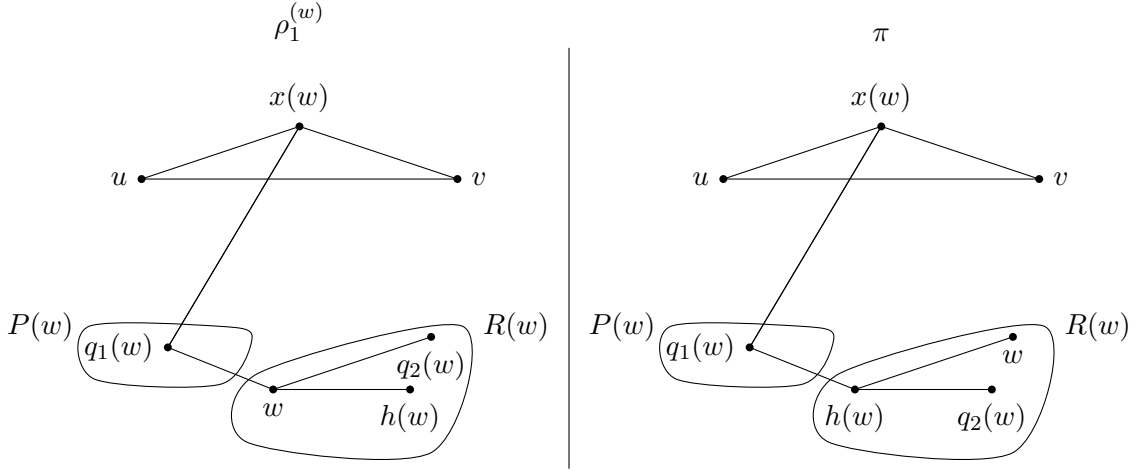

    Notice that if $\rho_1^{(w)} = \rho_1^{(w^\star)},$ then $h(w^\star) = h(w),$ $q_1(w^\star) = q_1(w),$ and $q_2(w^\star) = q_2(w)$ or $q_1(w^\star) = h(w),$ $q_2(w^\star) = q_1(w),$ and $h(w^\star) = q_2(w)$ or $q_2(w^\star) = h(w),$ $h(w^\star) = q_1(w),$ and $q_1(w^\star) = q_2(w),$ i.e., the positions of $h(w^\star), q_1(w^\star),$ and $q_2(w^\star)$ are one of three rotations of $h(w), q_1(w),$ and $q_2(w).$ 
    
    Similarly if $\pi = \rho_1^{(w^\star)},$ then $h(w^\star) = h(w),$ $q_1(w^\star) = q_1(w),$ and $q_2(w^\star) = w$ or $q_1(w^\star) = h(w),$ $q_2(w^\star) = q_1(w),$ and $h(w^\star) = w$ or $q_1(w^\star) = h(w),$ $q_2(w^\star) = q_1(w),$ and $h(w^\star) = w$ or $q_2(w^\star) = h(w),$ $h(w^\star) = q_1(w),$ and $q_1(w^\star) = w,$ i.e., the positions of $h(w^\star), q_1(w^\star),$ and $q_2(w^\star)$ are one of three rotations of $h(w), q_1(w),$ and $w$.
    
    Notice that when $\rho_1^{(w)} = \rho_1^{(w^\star)},$ then $w^\star$ is always adjacent to either $q_1(w)$ or $q_2(w)$ and $w^{\star \prime} = \sigma^{-1}(w^\star)$ is always adjacent to $h'(w)$, $q_1'(w)$, and $q_2'(w).$ Without loss of generality, let $w^\star$ be neighbors with $q_1(w).$ We can then use the following path to swap $u$ and $v.$
    \begin{gather*}
    (q_1'(w), w^{\star\prime}), (w^{\star\prime},h'(w)), (w',h'(w)), (h'(w),x'(w)), (x'(w),u'), (u',v'), \\(v',x'(w)), (h'(w),x'(w)), (w',h'(w)),(w^{\star\prime},h'(w)), (q_1'(w), w^{\star\prime}).
    \end{gather*}

    If $\pi = \rho_1^{(w^\star)}$ on the other hand, then we know that $w^\star$ is adjacent to $q_1(w)$ or $w^\star$ is adjacent to both $w$ and $h(w).$ In the first case, we can use the same path as above to swap $u$ and $v$.
    Otherwise, we can use the path
    \begin{gather*}
    (w^{\star \prime},h'(w)), (w',h'(w)), (h'(w),x'(w)), (x'(w),u'), (u',v'), \\(v',x'(w)), (h'(w),x'(w)), (w',h'(w)), (w^{\star \prime},h'(w)).
    \end{gather*}
    Notice that for different $w^\star_1$ and $w^\star_2$ the corresponding paths are disjoint because in the first path, all intermediate permutations satisfy that $w^\star_1$ changes position while $w^\star_2$ remains fixed and in the second path, all intermediate permutations satisfy that $w^\star_1$ remains fixed while $w^\star_2$ changes position. Therefore, if there are more than $k-1$ such $w^\star,$ there exist $k$ disjoint paths between $\sigma$ and $\sigma\circ (u,v)$.

    Thus for any $w,$ there are at most $k+3$ elements $w^\star$ including $w, q_1(w), q_2(w),$ and $h(w)$ which satisfy that some vertex in $\sigma(P(w))$ has more than one neighbor in $\sigma(R(w))$ and the path of $w^\star$ is not disjoint from the path of $w.$
    It follows that there are at most $(k-1)(k+3)$ elements satisfying that some vertex in $\sigma(P(w^\star))$ has more than one neighbor in $\sigma(R(w^\star))$ as any more would create $k$ disjoint paths.
\end{proof}
Combining the previous three lemmas, because $\lvert W \rvert > (k-1)(k+6),$ we know that there must exist a vertex, $w \in W$ such that $x'(w)$ is not adjacent to $w',$ $R(w)$ and $N[x'(w)]$ are disjoint, and no vertices in $\sigma(P(w))$ have more than one neighbor in $\sigma(R(w)).$ From this point, the same idea as the end of Bangachev's proof finishes this proof \cite[p.~18-19]{bangachev2022asymmetric}.
\subsection{\cref{kban2}}
As a result of \cref{starkconnectivity} we have \cref{kban2}, which we restate below for convenience. The proof is analogous to the one provided by Bangachev for \cref{ban2} \cite[p.~11-12]{bangachev2022asymmetric}. We use \cref{exchangeability} in place of \cite[Lemma 2.11]{bangachev2022asymmetric} and we invoke \cref{genwilson} together with the lower bound of $k$ on the minimum degree of $Y\lvert_Q$ as defined in \cite{bangachev2022asymmetric}.

\begin{customthm}{1.7}
    Let $X$ and $Y$ be two connected graphs on $n$ vertices satisfying
    \begin{itemize}
        \item $\delta(X) + \delta(Y) \geq n+k-1,$
        \item $\min\{\delta(X), \delta(Y)\}+2\max\{\delta(X),\delta(Y)\} \geq 2n,$
    \end{itemize}
    where $k \geq 2.$ Then $\FS(X,Y)$ is $k$-connected.
\end{customthm}

\section{Random Graphs} \label{section:random}
We now extend the results of Alon, Defant, and Kravitz \cite{alon2021extremal}, Milojevi\' c \cite{milojevic2022connectivity}, and Wang and Chen \cite{wang2023connectivity} where $X$ and $Y$ are random graphs $\mathcal{G}(n,p)$ by identifying a threshold probability $p$ above which $\FS(X,Y)$ is $k$-connected.
\subsection{Starcle Connectivity}
    We start by considering starcle graphs. Finding the connectivity of the friends-and-strangers graph formed by a starcle and a star graph is important as we will show that $X$ contains a starcle graph and $Y$ contains a star graph with high probability where $X\sim \mathcal{G}(n,p_1)$ and $Y\sim \mathcal{G}(n,p_2)$ which will then give a highly connected subgraph of $\FS(X,Y)$. We define the \textit{starcle graph} on $n$ vertices and diagonal tuple $(x_1, \ldots, x_{k-1})$ where $1 < x_1 < x_2 < \cdots < x_{k-1} < n-1$ and $k \leq n-2$ to be the graph with vertex set $[n]$ and edge set $\{(i,i+1):1\leq i<n\} \cup \{(1,n)\} \cup \{(x_i,n):1\leq i<k\}.$ Intuitively these graphs look like a cycle combined with a star where the center of the star lies at one of the vertices of the cycle. Let us prove the following lemma on the existence of disjoint paths in $\FS(X,\Star_n)$ where $X$ is a starcle graph.
    \begin{lemma}\label{starcle}
        Let $n\geq 3$ and $k\geq 2.$ Let $X$ be the starcle graph on $n$ vertices and diagonal tuple $(x_1, \ldots, x_{k-1}).$ Let $x_0 = 1$ and $x_k = n-1.$ Assume that there exists $1\leq i\leq k$ such that $x_i-x_{i-1}$ is odd, i.e., $X$ has an odd cycle. Also assume that $x_{i+1}-x_i\geq k$ for $0 \leq i \leq k-1.$ Let $G = \FS(X,\Star_n).$ Then there exist $k$ disjoint paths between any two distinct vertices $\sigma$ and $\rho$ where $\sigma^{-1}(n) = \rho^{-1}(n) = n.$ 
    \end{lemma}
    \begin{proof}
        Define $\pi_i\in\mathfrak{S}_{n-1}$ to be the cycle $(x_{i+1},x_{i+1}-1,\ldots,x_i)$ for $0 \leq i \leq k-1.$ We will show that these $\pi_i$ generate all permutations in $\mathfrak{S}_{n-1}.$ First, note that $G$ is connected by \cref{wilson} as $X$ is biconnected, not bipartite, not isomorphic to a cycle graph, and not isomorphic to $\theta_0.$ 
        
        Let $\sigma$ and $\rho$ be defined as above. There must exist a path between $\sigma$ and $\rho$ as $G$ is connected. Define $\tau_1, \tau_2, \ldots \tau_\ell$ to be a consecutive subset of the vertices of the path between $\sigma$ and $\rho$ such that $\tau_1^{-1}(n) = \tau_\ell^{-1}(n) = n$ and no $1 < i < \ell$ satisfies $\tau_i^{-1}(n) = n.$ Notice that $\tau_2^{-1}(n) = x_p$ and $\tau_{\ell-1}^{-1}(n)=x_q$ for some $p$ and $q.$ Assume without loss of generality that $p < q.$ Then $\tau_{1} = \tau_{\ell} \circ \pi_{q-1} \circ \pi_{q-2} \circ \cdots \circ \pi_p$, viewing the $\pi_i$ as elements in $\mathfrak{S}_{n-1}$ by adding $n$ as a fixed point. Therefore, the path from $\sigma$ to $\rho$ can be broken up into a composition of permutations of the form $\pi_i,$ and because $\sigma$ and $\rho$ are completely arbitrary aside from the condition that $\sigma^{-1}(n) = \rho^{-1}(n) = n$ we know $\pi_i$ must generate $\mathfrak{S}_{n-1}.$

        Next, note that the set of permutations $P = \{\pi_i : 0\leq i \leq k-1\}$ is a minimal generating set of $\mathfrak{S}_{n-1}.$
        Let $H$ be the Cayley graph with vertex set $\mathfrak{S}_{n-1}$ and edge set such that there exists an edge between $\tau$ and $\tau'$ if $\tau^{-1} \circ \tau'$ or $(\tau')^{-1}\circ\tau$ is in $P.$
        As shown by Godsil \cite{godsil1981cayley}, we then have that $H$ is $\lvert P \cup P^{-1}\rvert$-connected, where $P^{-1} = \{\pi_i^{-1} : 0 \leq i \leq k-1\}.$
        Because $x_{i+1}-x_i \geq 2,$ we know $\lvert P \cup P^{-1}\rvert = 2k,$ so $H$ is $2k$-connected.
        
        Let $J$ be the directed graph with vertex set $\mathfrak{S}_{n-1}$ and edge set such that there exists a directed edge from $\tau$ to $\tau'$ if $\tau^{-1} \circ \tau' \in P.$ Let us show that there exist $k$ disjoint directed paths from $\sigma\lvert_{[n-1]}$ to $\rho\lvert_{[n-1]}$ in $J.$ 
        
        Because Menger's theorem holds for directed graphs (e.g. see Göring \cite{goring2000digraphmenger}), the existence of $k$ disjoint paths is equivalent to showing that all vertex cuts between $\sigma\lvert_{[n-1]}$ to $\rho\lvert_{[n-1]}$ have at least $k$ vertices. Let $A$ be a subset of $V(H).$ Let $N_{\out}(A)$ be the set of vertices $\tau \not \in A$ such that there exists a directed edge in $J$ from a vertex in $A$ to $\tau.$ Let $N_{\xin}(A)$ be the set of vertices $\tau \not \in A$ such that there exists a directed edge in $J$ from $\tau$ to a vertex in $A.$ Showing that all vertex cuts have at least $k$ vertices is equivalent to showing that for all $A$ with $A \cup N_{\out}(A) \neq V(H),$ we have $\lvert N_\out(A) \rvert \geq k.$ If there is a cut of $k$ vertices, then any connected component of the remaining graph has $\lvert N_{\out}(A)\rvert < k,$ and any other connected component is not in its neighborhood, and if some set $A$ satisfies $\lvert N_{\out}(A)\rvert < k,$ then removing its neighborhood disconnects the graph.

        Let us show that $\lvert N_{\out}(A) \rvert = \lvert N_{\xin}(A) \rvert.$ If there were no directed edges between two vertices of $A,$ because for any vertex $\tau$ in $A,$ $\lvert N_{\xin}(\tau) \rvert = \lvert N_{\out}(\tau) \rvert = k,$ we have $\lvert N_{\out}(A) \rvert = \lvert N_{\xin}(A) \rvert.$ For every additional directed edge between two vertices of $A$, it removes an edge moving out of a vertex of $A$ and it removes a vertex moving into a vertex of $A.$
        Therefore with every additional directed edge, both $\lvert N_{\out}(A) \rvert$ and $\lvert N_{\xin}(A) \rvert$ decrease by $1.$ 
        
        Let $N(A) = N_{\out}(A) \cup N_{\xin}(A).$ Note that it suffices to show $\lvert N(A) \rvert \geq 2k$, as this implies $\lvert N_\out(A) \rvert \geq k.$ If $A \cup N(A) \neq V(H),$ then because $H$ is $2k$-connected, we know by definition that $\lvert N(A) \rvert \geq 2k.$ 
        
        We now consider when $A \cup N(A) = V(H).$ Let us show that the girth of $J$ is greater than or equal to $k.$ We define girth for a directed graph to be the minimum length of a directed cycle. We will do this by considering any cycle. Let the intermediate permutations of this cycle be $\tau_0, \tau_1, \ldots, \tau_{w}$ where $\tau_0 = \tau_w.$ Start by performing any permutation $\pi_i = (x_{i+1},\ldots,x_i)$ to move from $\tau_0$ to $\tau_1.$ There must be some vertex $v \in [n-1]$ such that $\tau_0^{-1}(v) \neq \tau_1^{-1}(v).$ After the series of permutations, person $v$ must be restored to their original position. 
        
        Notice that there must be some $\tau_j$ such that $\tau_j^{-1}(v) > \tau_{j+1}^{-1}(v)$ as otherwise, the position of $v$ will just keep increasing. Notice that the first $j$ where this is possible is $x_{i+1}-\tau_0^{-1}(v).$ If $\tau_{j+1}^{-1}(v) \geq x_{i+1},$ then $v$ must have reached at least position $x_{i+2}$ which requires at least $k$ applications of a permutation. Therefore, we only need to consider when $\tau_{j+1}^{-1}(v) \leq x_i.$ From this position, we can either increment the position of $v$ by $1$ or move $v$ back. Therefore, the minimum number of moves required to return $v$ to its original position is $\tau_0^{-1}(v) - x_i.$ Therefore, the minimum number of swaps is $x_{i+1}-x_i \geq k.$

        Recall we have $A\cup N_{\out}(A)\neq V(H)$ and $A\cup N(A)=V(H)$.
        Let $S=N_{\xin}(A) \setminus (N_{\xin}(A) \cap N_{\out}(A)).$ 
        Notice that we may assume $\lvert S \rvert < k$, as otherwise $\lvert N_{\out}(A) \setminus (N_{\xin}(A) \cap N_{\out}(A)) \rvert = \lvert S\rvert \geq k$ and thus $\lvert N(A)\rvert \geq 2k$, which would complete the proof.
        Let $\tau$ be a vertex in $S.$ This vertex must exist by the assumptions that $A \cup N_\out(A) \neq V(H)$ and $A\cup N(A)=V(H)$.
        Consider $N_\xin(\tau),$ the set of $k$ vertices where there exists a directed edge from that vertex to $\tau.$ Because $N_\xin(\tau) \cap A$ is empty by definition, we know that $N_\xin(\tau) \subseteq N(A).$ If $N_\xin(\tau)$ has no vertices in $S,$ then $N_\xin(\tau) \subseteq N_\out(A),$ so $\lvert N_\out(A)\rvert \geq k.$
        Otherwise, for every $\tau \in S,$ there must exist a directed edge to $\tau$ from another vertex in $S.$ Notice that this must create a cycle within $S$, but because $\lvert S \rvert < k,$ this is impossible as the girth of $J$ is greater than or equal to $k.$ This shows that for all $A$ with $A \cup N_\out(A) \neq V(H),$ we have $\lvert N_\out(A) \rvert \geq k.$
        Thus, there exist $k$ disjoint directed paths in $J$ from $\sigma\lvert_{[n-1]}$ to $\rho\lvert_{[n-1]}.$
 
        Notice that $\pi_i$ can be expressed as the following series of transpositions:
        $$(n,x_{i+1}) \circ (x_{i+1},x_{i+1}-1) \circ \cdots \circ (x_i+1,x_i) \circ (x_i, n).$$
        Notice that replacing each directed edge in these $k$ paths with friendly transpositions corresponding to the permutation in $P$ creates $k$ paths from $\sigma$ to $\rho.$ Note that these $k$ paths from $\sigma$ to $\rho$ are on a transposition level in $\FS(X,Y)$ as opposed to a permutation level in $H.$
        
        We will now show that these paths are disjoint in $\FS(X,Y)$. Recall that the type of a vertex $\sigma \in G$ is $\sigma^{-1}(n).$ If there is a vertex $\tau$ shared by two paths, then notice that the type of $\tau$ cannot be $n$ as then the $k$ directed paths in $J$ wouldn't be disjoint. Say the type of $\tau$ is $x_i$ for some $1 \leq i \leq k-1.$ Then in the two paths, $Q$ and $R,$ that share $\tau,$ the permutation in $Q$ just before or after $\tau$ must be of type $n.$ Similarly, the permutation just before or after $\tau$ in $R$ must be of type $n.$ These two vertices must be the same as they are created by the same transposition. Then, the $k$ directed paths in $J$ aren't disjoint. Finally, if $\tau$ is of any other type, the type determines the overall permutation $\pi_i$ that is being performed. Then by completing the permutation, both paths must end on the same vertex of type $n$ which contradicts that the $k$ directed paths in $J$ are disjoint. This completes the proof of the existence of $k$ disjoint paths.
    \end{proof}

We also consider when multiple starcle graphs are put together to form the $X$ graph in $\FS(X,\Star_n).$

\begin{lemma}\label{multifs}
    Let $X$ be a graph on $n\geq 4$ vertices. Assume $X$ has two subgraphs $A$ and $B$ such that $V(A) \cup V(B) = V(X)$. Let $G_1 = \FS(A,\Star_n\lvert_{V(A)})$ and $G_2 = \FS(B,\Star_n\lvert_{V(B)}).$ Assume $n \in V(X)$ is in both $V(A)$ and $V(B)$ and assume there exist $k\geq 2$ disjoint paths between any two vertices $\tau_1$ and $\tau_2$ in $G_1$ and similarly for $G_2$ where $\tau_1^{-1}(n) = \tau_2^{-1}(n) = n.$ Let $I = V(A) \cap V(B) \setminus \{n\}$ and suppose $\lvert I \rvert \geq k.$ Then, there exist $k$ disjoint paths between any two vertices $\sigma$ and $\rho$ in $G = \FS(X,\Star_n)$ where $\sigma^{-1}(n) = \rho^{-1}(n) = n.$
\end{lemma}
\begin{proof}
    We start by splitting the path from $\sigma$ to $\rho$ into two phases divided into smaller steps.
    Let the subgraph list of a person $i$ be the subset of $\{A,B\}$ consisting of the subgraphs that $\rho^{-1}(i)$ is in.
    The first phase is to move the people around such that they are in one of the subgraphs of their subgraph list. The second phase is to order the people within these subgraphs until they are in the right position.

    Let us start with the first phase. Our goal is to break up this phase into smaller steps where in each step, we are permuting the vertices of just one of the subgraphs. To achieve the first phase, it suffices to get all people whose subgraph lists only contain $B$ to be on vertices in $V(B)\setminus I.$
    We do this using the following algorithm.
    \begin{enumerate}
        \item Let $S$ be the set of people standing on vertices in $A$ whose subgraph list only contains $B.$
        \item Move $\ell = \min \{\lvert S \lvert, \lvert I \lvert)$ of these people onto positions in $I$ through a permutation in subgraph $A,$ maintaining the position of person $n.$ Notice that $n$ is in subgraph $A,$ so this permutation must be possible.
        \item Move the $\ell$ vertices from $I$ to positions in $V(B) \setminus I$ without adding any other people whose subgraph list only contains $B$ into $I$ through a permutation in subgraph $B$ maintaining the position of person $n.$ We must always be able to do this as the total number of people whose subgraph list consists of $B$ is $\lvert V(B) \setminus I \rvert.$ Notice that $n$ is in subgraph $B,$ so this permutation must be possible.
        \item Repeat from (1).
    \end{enumerate}

    For the second phase, we know all people standing on vertices in subgraph $A$ have a subgraph list that contains $A.$ Therefore, we can perform a permutation on subgraph $A$ to put all of them in their correct position. Finally we can perform a permutation in subgraph $B$ to put the remaining people in the right positions.

    Let the sequence of these permutations be $\pi_1, \pi_2, \ldots, \pi_r$ for some $r$ where $\pi_i \in V(G).$ Note that in all of the intermediate permutations except for $\sigma$ and $\rho$ in this path, the ordering of the people standing on vertices in $I$ does not matter.
    Because $\lvert I \rvert \geq k,$ there must exist at least $k$ distinct permutations of the vertices on $I.$
    We refer to any one of these $k$ permutations of $\pi_i$ as $\pi_{i,j}$ for $1\leq j \leq k.$

    Let us show that in $G_1$ and $G_2,$ from any permutation $\tau,$ for permutations $\tau_1, \tau_2, \ldots, \tau_k$ where $\tau,\tau_1,\ldots,\tau_k$ are all distinct, there exist disjoint paths from $\tau$ to $\tau_i$ where $\tau^{-1}(n) = \tau_i^{-1}(n)= n$ for all $1 \leq i \leq k.$
    Without loss of generality, we just consider $G_1.$ Add a vertex $\tau'$ to form a new graph $G_1^\star$ where $\tau'$ is adjacent to $\tau_i$ for all $1 \leq i \leq k.$ 
    
    Let us show there exist $k$ disjoint paths from $\tau$ to $\tau'.$ If we consider any cut of size $k-1,$ one of the $\tau_i$ cannot have been removed, say $i = j.$ Then there exist $k$ disjoint paths from $\tau$ to $\tau_j,$ so removing $k-1$ vertices cannot disconnect them. As a result, there exists a path from $\tau$ to $\tau',$ so by Menger's Theorem \cite{diestelblocks}, there exist $k$ disjoint paths from $\tau$ to $\tau'$ which means there exist disjoint paths from $\tau$ to $\tau_i$ for $1\leq i \leq k.$

    As a result, we have that there exist disjoint paths from $\pi_{i,j}$ and $\pi_{i+1,\ell}$ for fixed $1 \leq i \leq r-1$ and fixed $1\leq j \leq k$ and variable $1\leq \ell \leq k.$ There also exist disjoint paths from $\sigma$ to $\pi_{1,j}$ for $1\leq j \leq k.$

    Consider a vertex cut of $k-1$ vertices. Because $k-1$ are removed, there must exist a path from $\sigma$ to $\pi_{1,j_1}$ for some $j_1.$ Similarly, there must exist a path from $\pi_{1,j_1}$ to $\pi_{2,j_2}$ for some $j_2.$ We can continue this process until we reach $\pi_{r,j_r}$ for some $j_r$. Notice that there exist $k$ disjoint paths between $\pi_{r,j_r}$ and $\rho,$ so there must exist some path between them. Therefore, no cut of $k-1$ vertices disconnects $\sigma$ and $\rho,$ so there must exist $k$ disjoint paths between the two.
\end{proof}

\subsection{Main Proof}
We now proceed with the main proof of \cref{kconnectivitywhp}, which we restate below for convenience.
\begin{customthm}{1.9}
    Fix some small $\varepsilon > 0$ and a positive integer $k>1.$ Let $X$ and $Y$ be independently-chosen random graphs in $\mathcal{G}(n,p_1)$ and $\mathcal{G}(n,p_2)$ where $p_1=p_1(n)$ and $p_2=p_2(n)$ both depend on $n.$ If
    $$p_1p_2\leq \varepsilon\frac{\log n}{n},$$
    then $\FS(X,Y)$ is disconnected with high probability. If $p_1p_2 \geq p_0^2$ and $p_1,p_2\geq p_0/\ell$ where
    $$p_0\geq \frac{\exp((k+7)/4\cdot(\log n)^{2/3})}{n^{1/2}},$$
    and $\ell = \frac{(\log n)^{2/3}}{2}$ then $\FS(X,Y)$ is $k$-connected with high probability.
\end{customthm}

We first show the disconnected with high probability condition and then show the connected with high probability condition.
\subsubsection{Disconnected with high probability}
We extend \cref{milojevicdisconnectedwhp} for an asymmetric condition.

\begin{proposition}
    There exists a constant $\varepsilon > 0$ with the following property. For a large positive integer $n$ and $p_1$ and $p_2$ such that 
    $$p_1 p_2 \leq \frac{\varepsilon\log n}{n} ,$$
    if we choose $X$ from $\mathcal{G}(n,p_1)$ and $Y$ from $\mathcal{G}(n,p_2),$ then $\FS(X,Y)$ has an isolated vertex with high probability.
\end{proposition}
\begin{proof}
    As mentioned by Alon, Defant, and Kravitz \cite{alon2021extremal}, an isolated vertex in $\FS(X,Y)$ corresponds to a packing of $X$ and $Y$ in $K_n$. Plugging in $k=2$ to Theorem $1$ of Bollob\' as, Janson, and Scott \cite{bollobas2016packing} yields the following statement: there exists a constant $\varepsilon > 0$ such that for $p_1$ and $p_2$ where $p_1 p_2 \leq \varepsilon \log n / n,$ there exists a packing of $X$ and $Y$ in $K_n$ with high probability. Therefore $\FS(X,Y)$ has an isolated vertex with high probability.
\end{proof}

\subsubsection{Connected with high probability}
Our goal is to prove the following proposition.

\begin{proposition} \label{wangchenkconnectedwhp}
    Let $X$ and $Y$ be independently-chosen random graphs in $\mathcal{G}(n,p_1)$ and $\mathcal{G}(n,p_2),$ respectively. Let $p_0 = \frac{\exp((k+7)/4\cdot (\log n)^{2/3})}{n^{1/2}}$ and $\ell = \frac{(\log n)^{2/3}}{2}.$ If
    $$p_1p_2 \geq p_0^2 \quad \text{and} \quad p_1,p_2\geq \frac{1}{\ell}p_0,$$
    then $\FS(X,Y)$ is $k$-connected with high probability.
\end{proposition}

We start by restating definitions established by Alon, Defant, and Kravitz \cite{alon2021extremal} and Wang and Chen \cite{wang2023connectivity}.

Let $m$ be a positive integer. Let $G$ and $H$ be graphs with vertex set $[m].$
Let $X$ and $Y$ be graphs on $n$ vertices, and let $\sigma \in V(\FS(X,Y))$ be a permutation. Let $V_1, V_2, \ldots, V_m$ be pairwise disjoint subsets of $Y.$ 

We say that $(G,H)$ is \textit{embeddable} in $(X,Y)$ with respect to sets $V_1, \ldots, V_m$ and permutation $\sigma$ if there exist $v_i \in V_i$ for all $i$ such that if $(i,j)\in E(H),$ then $(v_i,v_j) \in E(Y)$ and if $(i,j) \in E(G)$ then $(\sigma^{-1}(v_i),\sigma^{-1}(v_j)) \in E(X).$

Let $q_1,q_2,\ldots, q_m$ be nonnegative integers such that $q_1 + q_2 + \cdots + q_m \leq n.$ We say $(G,H)$ is \textit{$(q_1,q_2,\ldots, q_m)$-embeddable} if for every list $V_1, V_2, \ldots, V_m$ where $V_i$ contains $q_i$ vertices and every permutation $\sigma$, $(G,H)$ is embeddable in $(X,Y)$ with respect to sets $V_1, \ldots, V_m$ and permutation $\sigma.$

We start with the following lemma of Wang and Chen \cite{wang2023connectivity} which is just an asymmetric version of Lemma 3.2 of Alon, Defant, and Kravitz \cite{alon2021extremal}.

\begin{lemma} [\protect{\cite[Lemma 3.3]{wang2023connectivity}}] \label{embeddablelemma}
    Let $m,n,q_1,\ldots,q_m$ be positive integers such that $Q:=q_1+\cdots+q_m\leq n,$ and let $G$ and $H$ be two graphs on the vertex set $[m].$ Let $X$ and $Y$ be independently chosen random graphs in $\mathcal{G}(n,p_1)$ and $\mathcal{G}(n,p_2)$ where $p_1$ and $p_2$ are functions of $n.$ If for every set $J\subseteq [m]$ satisfying $\lvert E(G\lvert_J)\rvert + \lvert E(H\lvert_J)\rvert\geq 1$ we have
    $$p_1^{\lvert E(G\lvert_J)\rvert}p_2^{\lvert E(H\lvert_J)\rvert}\prod_{j\in J} q_j \geq 3\cdot2^{m+1}Q\log n,$$
    then the probability that the pair $(G,H)$ is $(q_1,\ldots,q_m)$-embeddable in $\FS(X,Y)$ is at least $1-n^{-Q}.$
\end{lemma}

We now follow the strategy explained in Alon, Defant and Kravitz \cite{alon2021extremal} and Wang and Chen \cite{wang2023connectivity}. Let $n$ be a large integer, and let 
$$m=\left\lfloor (\log n)^{2/3} \right\rfloor.$$
We will construct specific graphs $G^\star$ and $H^\star$ on the vertex set $[m+2]$ such that there exist $k$ disjoint paths in $\FS(G^\star,H^\star)$ from the identity permutation to the permutation where $m+1$ and $m+2$ are swapped.
Then we will use \cref{embeddablelemma} in order to show that for graphs $X$ and $Y$ independently-chosen in $\mathcal{G}(n,p_1)$ and $\mathcal{G}(n,p_2),$ vertices $u, v \in V(Y),$ and permutation $\sigma \in V(\FS(X,Y))$, with high probability there exists an embedding $\varphi$ of $G^\star$ into $X$ and an embedding $\psi$ of $H^\star$ in $Y,$ such that $\varphi(m+1) = \sigma^{-1}(u),$ $\varphi(m+2) = \sigma^{-1}(v),$ and $\varphi\circ \varepsilon = \sigma \circ \psi$ where $\varepsilon$ is the identity map from $V(G^\star)$ to $V(H^\star).$ This implies that with high probability, there exist $k$-disjoint paths from $\sigma$ to $(u,v)\circ \sigma.$
The proof will then follow from \cref{exchangeability}.

We start by describing the graphs $G^\star$ and $H^\star.$ Note that $\ell$ is essentially $\lfloor \sqrt{m}/2 \rfloor$
Note that because floors do not matter for asymptotics, we omit them for the rest of the section. Let the elements of $[m]$ be
$$w,w_1,\ldots,w_{k-1},x_1,\ldots,x_\ell,y_1,\ldots,y_\ell,z_1,\ldots,z_{m-2\ell-k}.$$
Let $w = m$ without loss of generality (though we will still be using $w$ to refer to this vertex). Let us describe $H^\star$ first. Let it have edge set containing all $(w,v)$ for all other vertices $v$, all $(m+1,x_i)$ for $1 \leq i \leq \ell,$ and all $(m+2,y_i)$ for $1 \leq i \leq \ell.$ 
We also define $H^{\star\star} = H^\star \lvert_{[m]}.$
For $G^\star,$ we start constructing $G^{\star\star}$ by arranging the elements of $[m]$ in a cycle such that the vertices $z_1, \ldots, z_{12}$ appear in that order. We then add the edges $(z_1,z_6), (z_2,z_4),(z_7,z_{12}),(z_8,z_{10}),$ and $(w,w_i)$ for $1 \leq i \leq k-1.$ We make sure $G^\star$ has the following properties. Note that in the following points, distance refers to distance along the large cycle not including the additional edges.
\begin{itemize}
    \item $G^\star$ contains the edges $(z_4,z_5), (z_5,z_6), (z_{10},z_{11}),$ and $(z_{11},z_{12}).$ Thus $z_4,z_5,$ and $z_6$ occur consecutively and $z_{10},z_{11},$ and $z_{12}$ occur consecutively.
    \item The clockwise distance between $z_3$ and $z_5$ is $\ell-1.$ Similarly, the clockwise distance between $z_9$ and $z_{11}$ is $\ell-1.$
    \item The $2\ell+k$ vertices $w,w_1,\ldots,w_{k-1},x_1,\ldots,x_\ell,y_1,\ldots,y_\ell$ are placed such that the distance between any two of them and the distance between any one of them and one of $z_3, z_5,z_9,z_{11}$ is at least $m/(3\ell)$ which is at least $k$ for sufficiently large $n.$  
    \item The $k$ vertices $w,w_1,\ldots,w_{k-1}$ are placed in clockwise order around the cycle such that they lie between $z_1$ and $z_2$ or between $z_7$ and $z_8.$ The vertices $w$ and $w_1$ are between $z_7$ and $z_8$ in the clockwise direction, and the clockwise distance between $w$ and $w_1$ is an even number. The total number of vertices on the clockwise path from $z_1$ to $z_2$ or the clockwise path from $z_7$ to $z_8$ is at least $k+1.$
    \item For $k\geq 4,$ there exist at least two vertices $w_i$ and $w_j$ which are on the clockwise path from $z_1$ to $z_2.$ If $k = 2,$ there are no such vertices, and if $k = 3,$ there is one such vertex.
    \item The girth of the entire graph $G^{\star\star}$ is at least $m/(k+6),$ and $m$ is sufficiently large so that $m/(k+6)\geq k+2.$ 
\end{itemize}
The graph $G^\star$ is obtained from $G^{\star \star}$ by adding the edges $(m+1,z_3),(m+1,z_{11}),(m+2,z_5),(m+2,z_9),$ and $(m+1,m+2).$ From now on, we refer to the cycle in $G^{\star \star}$ containing all vertices in $[m]$ as the \textit{large cycle}.

\begin{figure}[htbp]
    \centering
    \vcenteredhbox{\begin{tikzpicture}[scale=3,
    dot/.style = {circle, fill, minimum size=4.5pt, 
              inner sep=0pt, outer sep=0pt}]
        \draw [red,ultra thick,domain=5:175] plot ({cos(\x)}, {sin(\x)});
        \draw [red,ultra thick,domain=185:355] plot ({cos(\x)}, {sin(\x)});
        \coordinate (z12) at ({cos(-5)}, {sin(-5)});
        \coordinate (z11) at ({cos(0)}, {sin(0)});
        \coordinate (z10) at ({cos(5)}, {sin(5)});
        \coordinate (z9) at ({cos(15)}, {sin(15)});
        \coordinate (z8) at ({cos(45)}, {sin(45)});
        \coordinate (z7) at ({cos(160)}, {sin(160)});
        \coordinate (z6) at ({cos(175)}, {sin(175)});
        \coordinate (z5) at ({cos(180)}, {sin(180)});
        \coordinate (z4) at ({cos(185)}, {sin(185)});
        \coordinate (z3) at ({cos(195)}, {sin(195)});
        \coordinate (z2) at ({cos(225)}, {sin(225)});
        \coordinate (z1) at ({cos(-20)}, {sin(-20)});
        \coordinate (w) at ({cos(90)}, {sin(90)});
        \coordinate (w1) at ({cos(50)}, {sin(50)});
        \coordinate (w2) at ({cos(-30)}, {sin(-30)});
        \coordinate (w3) at ({cos(-70)}, {sin(-70)});
        \coordinate (w4) at ({cos(-110)}, {sin(-110)});
        \coordinate (w5) at ({cos(130)}, {sin(130)});
        \coordinate (m1) at (-0.4,0);
        \coordinate (m2) at (0.4,0);
        \node [dot] at (z1) [label=below right: $z_1$] {};
        \node [dot] at (z2) [label=below left: $z_2$] {};
        \node [dot] at (z3) [label=below left: $z_3$] {};
        \node [dot] at (z4) [label={left,yshift=-0.1 cm}: $z_4$] {};
        \node [dot] at (z5) [label=left: $z_5$] {};
        \node [dot] at (z6) [label={left,yshift=0.1 cm}: $z_6$] {};
        \node [dot] at (z7) [label=above left: $z_7$] {};
        \node [dot] at (z8) [label=above right: $z_8$] {};
        \node [dot] at (z9) [label=above right: $z_9$] {};
        \node [dot] at (z10) [label={right,yshift=0.1 cm}: $z_{10}$] {};
        \node [dot] at (z11) [label=right: $z_{11}$] {};
        \node [dot] at (z12) [label={right,yshift=-0.1 cm}: $z_{12}$] {};
        \node [teal, dot] at (w) [label={above}: $w$] {};
        \node [teal, dot] at (w1) [label={above}: $w_{i_1}$] {};
        \node [teal, dot] at (w2) [label={below right}: $w_{i_2}$] {};
        \node [teal, dot] at (w3) [label={below}: $w_{i_3}$] {};
        \node [teal, dot] at (w4) [label={below}: $w_{i_4}$] {};
        \node [teal, dot] at (w5) [label={above left}: $w_{i_5}$] {};
        \node [blue, dot] at (m1) [label={above}: $m+1$] {};
        \node [blue, dot] at (m2) [label={below}: $m+2$] {};
        \draw[very thick] (z1)--(z6) (z12)--(z7) (z10)--(z8) (z4)--(z2) (z4)--(z5)--(z6) (z10)--(z11)--(z12);
        \draw [very thick, teal] (w)--(w1) (w)--(w2) (w)--(w3) (w)--(w4) (w)--(w5);
        \draw [very thick, blue] (m1)--(z3) (m2)--(z9) (m1)--(m2);
        \draw[very thick, blue] plot [smooth, tension=1] coordinates {(m2) (-0.4,0.3) (z5)};
        \draw[very thick, blue] plot [smooth, tension=1] coordinates {(m1) (0.4,-0.3) (z11)};
    \end{tikzpicture}}
    \hspace{0.8 cm}
    \vcenteredhbox{
    \begin{tikzpicture}[scale=2,
    dot/.style = {circle, fill, minimum size=4.5pt, 
              inner sep=0pt, outer sep=0pt}]
        \coordinate (w) at (0,0);
        \coordinate (wk) at ({cos(80)}, {sin(80)});
        \coordinate (w1) at ({cos(100)}, {sin(100)});
        \coordinate (y1) at ({cos(50)}, {sin(50)});
        \coordinate (y2) at ({cos(35)}, {sin(35)});
        \coordinate (yl) at ({cos(10)}, {sin(10)});
        \coordinate (x1) at ({cos(180-50)}, {sin(180-50)});
        \coordinate (x2) at ({cos(180-35)}, {sin(180-35)});
        \coordinate (xl) at ({cos(180-10)}, {sin(180-10)});
        \coordinate (z1) at ({cos(190)}, {sin(190)});
        \coordinate (z2) at ({cos(205}, {sin(205)});
        \coordinate (z3) at ({cos(220)}, {sin(220)});
        \coordinate (zm1) at ({cos(180-190)}, {sin(180-190)});
        \coordinate (zm2) at ({cos(180-205}, {sin(180-205)});
        \coordinate (zm3) at ({cos(180-220)}, {sin(180-220)});
        \coordinate (m1) at ({1.5*cos(180-35)}, {1.5*sin(180-35)});
        \coordinate (m2) at ({1.5*cos(35)}, {1.5*sin(35)});
        \node [dot] at (w) [label=below: $w$] {};
        \node [dot] at (w1) [label=above: $w_1$] {};
        \node [dot] at (wk) [label=above: $w_k$] {};
        \node [dot] at (x1) [label=above: $x_1$] {};
        \node [dot] at (x2) [label={left}: $x_2$] {};
        \node [dot] at (xl) [label=left: $x_l$] {};
        \node [dot] at (y1) [label=above: $y_1$] {};
        \node [dot] at (y2) [label=right: $y_2$] {};
        \node [dot] at (yl) [label=right: $y_l$] {};
        \node [dot] at (z1) [label=left: $z_1$] {};
        \node [dot] at (z2) [label=left: $z_2$] {};
        \node [dot] at (z3) [label=below: $z_3$] {};
        \node [dot] at (zm1) [label=right: $z_{m-2\ell-k}$] {};
        \node [dot] at (zm2) [label=right: $z_{m-2\ell-k-1}$] {};
        \node [dot] at (zm3) [label=right: $z_{m-2\ell-k-2}$] {};
        \node [blue, dot] at (m1) [label=above: $m+1$] {};
        \node [blue, dot] at (m2) [label=above: $m+2$] {};
        \draw [thick] (w)--(x1) (w)--(x2) (w)--(xl) (w)--(y1) (w)--(y2) (w)--(yl) (w)--(z1) (w)--(z2) (w)--(z3) (w)--(zm1) (w)--(zm2) (w)--(zm3) (w)--(w1) (w)--(wk);
        \draw [black, very thick, dotted,domain=220:320] plot ({cos(\x)}, {sin(\x)});
        \draw [black, very thick, dotted,domain=10:35] plot ({cos(\x)}, {sin(\x)});
        \draw [black, very thick, dotted,domain=180-35:180-10] plot ({cos(\x)}, {sin(\x)});
        \draw [black, very thick, dotted,domain=180-35:180-10] plot ({cos(\x)}, {sin(\x)});
        \draw [black, very thick, dotted,domain=80:100] plot ({cos(\x)}, {sin(\x)});
        \draw [blue, very thick] (m1)--(x1) (m1)--(x2) (m1)--(xl) (m2)--(y1) (m2)--(y2) (m2)--(yl);
    \end{tikzpicture}
    }
    \caption{Schematic diagrams $G^\star$ (left) and $H^\star$ (right). Blue, black, and green curves represent edges in the graph. Red curves represent collection of vertices and edges. The vertices $x_1,\ldots,x_\ell$ and $y_1,\ldots,y_\ell$ are not marked in $G^\star.$ The diagram of $G^\star$ shows $w$ and $w_i$ for $i \in \{i_1,\ldots i_5\}$ where $1\leq i_1 < \ldots < i_5 \leq k-1.$}
    \label{fig:gstarhstar}
\end{figure}
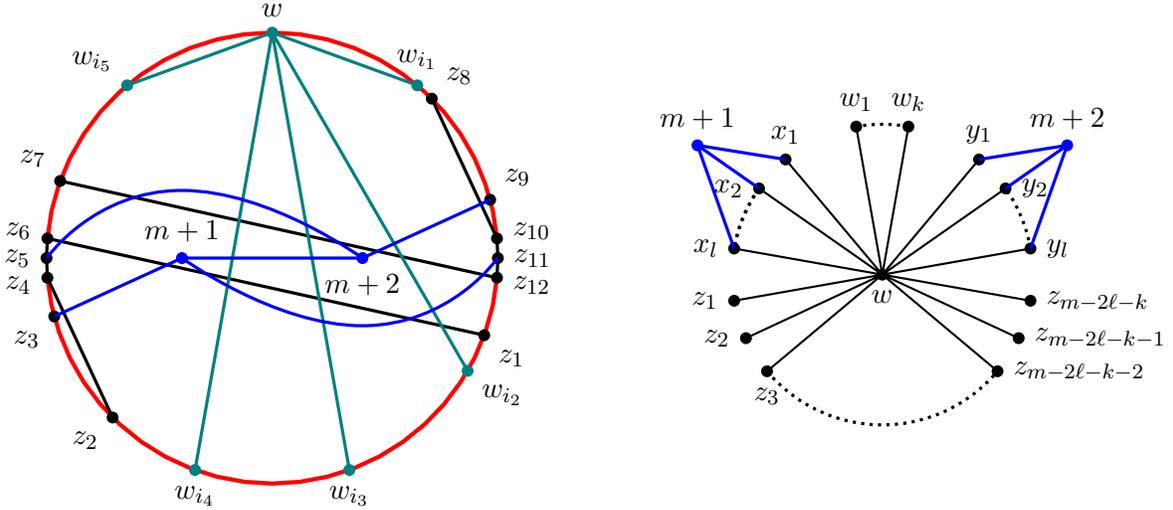
We now show that we can form $k$ disjoint paths that swap $m+1$ and $m+2.$
\begin{lemma}
    Let $G^\star$ and $H^\star$ be the graphs constructed above (see \cref{fig:gstarhstar}). Let $\sigma$ be the identity permutation. There exist $k$ disjoint paths from $\sigma$ to $\sigma \circ (m+1,m+2)$ in $\FS(G^\star,H^\star).$
\end{lemma}
\begin{proof}
    Let $G^{\star \star} = G^\star \lvert_{[m]}$ and $G^{\star \star \star} = G^\star \lvert_{[m+2]\setminus\{z_5,z_{11}\}}.$ Note that $G^{\star \star}$ contains a starcle graph on $m$ vertices with diagonal tuple $(w_1, \ldots, w_{k-1}).$ There is an odd cycle as the distance between $w$ and $w_1$ is even along the large cycle.
    By \cref{starcle}, we know that in $\FS(G^{\star \star}, H^{\star\star}),$ there exist $k$ disjoint paths between any two vertices $\rho_1$ and $\rho_2$ where $\rho_1^{-1}(w) = \rho_2^{-1}(w) = w.$

    As for $G^{\star \star \star},$ consider the graph $G'$ created by removing all edges of the form $(w,w_i)$ for $1 \leq i \leq k-1.$ Notice that $E(G')$ is the non-disjoint union of three cycles, depicted in \cref{fig:a'b'c'}, which we call $A'$, $B'$, and $C'$.

    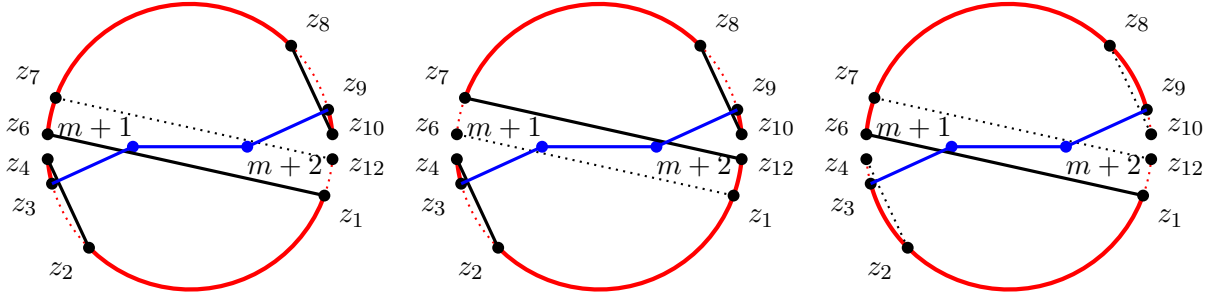
\begin{figure}[htbp]
    \centering
    \vcenteredhbox{\begin{tikzpicture}[scale=1.9,
    dot/.style = {circle, fill, minimum size=4.5pt, 
              inner sep=0pt, outer sep=0pt}]
        \draw [red,ultra thick,domain=5:15] plot ({cos(\x)}, {sin(\x)});
        \draw [red,dotted,thick,domain=15:45] plot ({cos(\x)}, {sin(\x)});
        \draw [red,ultra thick,domain=45:160] plot ({cos(\x)}, {sin(\x)});
        \draw [red,ultra thick,domain=160:175] plot ({cos(\x)}, {sin(\x)});
        \draw [red,ultra thick,domain=185:195] plot ({cos(\x)}, {sin(\x)});
        \draw [red,dotted,thick,domain=195:225] plot ({cos(\x)}, {sin(\x)});
        \draw [red,ultra thick,domain=225:340] plot ({cos(\x)}, {sin(\x)});
        \draw [red,dotted,thick,domain=340:355] plot ({cos(\x)}, {sin(\x)});
        \coordinate (z12) at ({cos(-5)}, {sin(-5)});
        \coordinate (z11) at ({cos(0)}, {sin(0)});
        \coordinate (z10) at ({cos(5)}, {sin(5)});
        \coordinate (z9) at ({cos(15)}, {sin(15)});
        \coordinate (z8) at ({cos(45)}, {sin(45)});
        \coordinate (z7) at ({cos(160)}, {sin(160)});
        \coordinate (z6) at ({cos(175)}, {sin(175)});
        \coordinate (z5) at ({cos(180)}, {sin(180)});
        \coordinate (z4) at ({cos(185)}, {sin(185)});
        \coordinate (z3) at ({cos(195)}, {sin(195)});
        \coordinate (z2) at ({cos(225)}, {sin(225)});
        \coordinate (z1) at ({cos(-20)}, {sin(-20)});
        \coordinate (w) at ({cos(90)}, {sin(90)});
        \coordinate (w1) at ({cos(50)}, {sin(50)});
        \coordinate (w2) at ({cos(-30)}, {sin(-30)});
        \coordinate (w3) at ({cos(-70)}, {sin(-70)});
        \coordinate (w4) at ({cos(-110)}, {sin(-110)});
        \coordinate (w5) at ({cos(130)}, {sin(130)});
        \coordinate (m1) at (-0.4,0);
        \coordinate (m2) at (0.4,0);
        \node [dot] at (z1) [label=below right: $z_1$] {};
        \node [dot] at (z2) [label=below left: $z_2$] {};
        \node [dot] at (z3) [label=below left: $z_3$] {};
        \node [dot] at (z4) [label={left,yshift=-0.1 cm}: $z_4$] {};
        \node [dot] at (z6) [label={left,yshift=0.1 cm}: $z_6$] {};
        \node [dot] at (z7) [label=above left: $z_7$] {};
        \node [dot] at (z8) [label=above right: $z_8$] {};
        \node [dot] at (z9) [label=above right: $z_9$] {};
        \node [dot] at (z10) [label={right,yshift=0.1 cm}: $z_{10}$] {};
        \node [dot] at (z12) [label={right,yshift=-0.1 cm}: $z_{12}$] {};
        \node [blue, dot] at (m1) [label={above,xshift=-0.5 cm,yshift=-0.1 cm}: $m+1$] {};
        \node [blue, dot] at (m2) [label={below,xshift=0.5 cm,yshift=0.1 cm}: $m+2$] {};
        \draw[very thick] (z1)--(z6) (z10)--(z8) (z4)--(z2);
        \draw[dotted, thick] (z12)--(z7);
        \draw [very thick, blue] (m1)--(z3) (m2)--(z9) (m1)--(m2);
    \end{tikzpicture}}
    \vcenteredhbox{\begin{tikzpicture}[scale=1.9,
    dot/.style = {circle, fill, minimum size=4.5pt, 
              inner sep=0pt, outer sep=0pt}]
        \draw [red,ultra thick,domain=5:15] plot ({cos(\x)}, {sin(\x)});
        \draw [red,dotted,thick,domain=15:45] plot ({cos(\x)}, {sin(\x)});
        \draw [red,ultra thick,domain=45:160] plot ({cos(\x)}, {sin(\x)});
        \draw [red,dotted, thick,domain=160:175] plot ({cos(\x)}, {sin(\x)});
        \draw [red, ultra thick,domain=185:195] plot ({cos(\x)}, {sin(\x)});
        \draw [red, dotted,thick,domain=195:225] plot ({cos(\x)}, {sin(\x)});
        \draw [red,ultra thick,domain=225:340] plot ({cos(\x)}, {sin(\x)});
        \draw [red,ultra thick,domain=340:355] plot ({cos(\x)}, {sin(\x)});
        \coordinate (z12) at ({cos(-5)}, {sin(-5)});
        \coordinate (z11) at ({cos(0)}, {sin(0)});
        \coordinate (z10) at ({cos(5)}, {sin(5)});
        \coordinate (z9) at ({cos(15)}, {sin(15)});
        \coordinate (z8) at ({cos(45)}, {sin(45)});
        \coordinate (z7) at ({cos(160)}, {sin(160)});
        \coordinate (z6) at ({cos(175)}, {sin(175)});
        \coordinate (z5) at ({cos(180)}, {sin(180)});
        \coordinate (z4) at ({cos(185)}, {sin(185)});
        \coordinate (z3) at ({cos(195)}, {sin(195)});
        \coordinate (z2) at ({cos(225)}, {sin(225)});
        \coordinate (z1) at ({cos(-20)}, {sin(-20)});
        \coordinate (w) at ({cos(90)}, {sin(90)});
        \coordinate (w1) at ({cos(50)}, {sin(50)});
        \coordinate (w2) at ({cos(-30)}, {sin(-30)});
        \coordinate (w3) at ({cos(-70)}, {sin(-70)});
        \coordinate (w4) at ({cos(-110)}, {sin(-110)});
        \coordinate (w5) at ({cos(130)}, {sin(130)});
        \coordinate (m1) at (-0.4,0);
        \coordinate (m2) at (0.4,0);
        \node [dot] at (z1) [label=below right: $z_1$] {};
        \node [dot] at (z2) [label=below left: $z_2$] {};
        \node [dot] at (z3) [label=below left: $z_3$] {};
        \node [dot] at (z4) [label={left,yshift=-0.1 cm}: $z_4$] {};
        \node [dot] at (z6) [label={left,yshift=0.1 cm}: $z_6$] {};
        \node [dot] at (z7) [label=above left: $z_7$] {};
        \node [dot] at (z8) [label=above right: $z_8$] {};
        \node [dot] at (z9) [label=above right: $z_9$] {};
        \node [dot] at (z10) [label={right,yshift=0.1 cm}: $z_{10}$] {};
        \node [dot] at (z12) [label={right,yshift=-0.1 cm}: $z_{12}$] {};
        \node [blue, dot] at (m1) [label={above,xshift=-0.5 cm,yshift=-0.1 cm}: $m+1$] {};
        \node [blue, dot] at (m2) [label={below,xshift=0.5 cm,yshift=0.1 cm}: $m+2$] {};
        \draw[very thick] (z12)--(z7) (z10)--(z8) (z4)--(z2);
        \draw[dotted,thick] (z6)--(z1);
        \draw [very thick, blue] (m1)--(z3) (m2)--(z9) (m1)--(m2);
    \end{tikzpicture}}
    \vcenteredhbox{\begin{tikzpicture}[scale=1.9,
    dot/.style = {circle, fill, minimum size=4.5pt, 
              inner sep=0pt, outer sep=0pt}]
        \draw [red,dotted, thick,domain=5:15] plot ({cos(\x)}, {sin(\x)});
        \draw [red,ultra thick,domain=15:45] plot ({cos(\x)}, {sin(\x)});
        \draw [red,ultra thick,domain=45:160] plot ({cos(\x)}, {sin(\x)});
        \draw [red,ultra thick,domain=160:175] plot ({cos(\x)}, {sin(\x)});
        \draw [red,dotted, thick,domain=185:195] plot ({cos(\x)}, {sin(\x)});
        \draw [red,ultra thick,domain=195:225] plot ({cos(\x)}, {sin(\x)});
        \draw [red,ultra thick,domain=225:340] plot ({cos(\x)}, {sin(\x)});
        \draw [red,dotted, thick,domain=340:355] plot ({cos(\x)}, {sin(\x)});
        \coordinate (z12) at ({cos(-5)}, {sin(-5)});
        \coordinate (z11) at ({cos(0)}, {sin(0)});
        \coordinate (z10) at ({cos(5)}, {sin(5)});
        \coordinate (z9) at ({cos(15)}, {sin(15)});
        \coordinate (z8) at ({cos(45)}, {sin(45)});
        \coordinate (z7) at ({cos(160)}, {sin(160)});
        \coordinate (z6) at ({cos(175)}, {sin(175)});
        \coordinate (z5) at ({cos(180)}, {sin(180)});
        \coordinate (z4) at ({cos(185)}, {sin(185)});
        \coordinate (z3) at ({cos(195)}, {sin(195)});
        \coordinate (z2) at ({cos(225)}, {sin(225)});
        \coordinate (z1) at ({cos(-20)}, {sin(-20)});
        \coordinate (w) at ({cos(90)}, {sin(90)});
        \coordinate (w1) at ({cos(50)}, {sin(50)});
        \coordinate (w2) at ({cos(-30)}, {sin(-30)});
        \coordinate (w3) at ({cos(-70)}, {sin(-70)});
        \coordinate (w4) at ({cos(-110)}, {sin(-110)});
        \coordinate (w5) at ({cos(130)}, {sin(130)});
        \coordinate (m1) at (-0.4,0);
        \coordinate (m2) at (0.4,0);
        \node [dot] at (z1) [label=below right: $z_1$] {};
        \node [dot] at (z2) [label=below left: $z_2$] {};
        \node [dot] at (z3) [label=below left: $z_3$] {};
        \node [dot] at (z4) [label={left,yshift=-0.1 cm}: $z_4$] {};
        \node [dot] at (z6) [label={left,yshift=0.1 cm}: $z_6$] {};
        \node [dot] at (z7) [label=above left: $z_7$] {};
        \node [dot] at (z8) [label=above right: $z_8$] {};
        \node [dot] at (z9) [label=above right: $z_9$] {};
        \node [dot] at (z10) [label={right,yshift=0.1 cm}: $z_{10}$] {};
        \node [dot] at (z12) [label={right,yshift=-0.1 cm}: $z_{12}$] {};
        \node [blue, dot] at (m1) [label={above,xshift=-0.5 cm,yshift=-0.1 cm}: $m+1$] {};
        \node [blue, dot] at (m2) [label={below,xshift=0.5 cm,yshift=0.1 cm}: $m+2$] {};
        \draw[very thick] (z1)--(z6);
        \draw[dotted,thick] (z12)--(z7) (z10)--(z8) (z4)--(z2);
        \draw [very thick, blue] (m1)--(z3) (m2)--(z9) (m1)--(m2);
    \end{tikzpicture}}
    \caption{Cycles $A', B', C'$ (in left to right order) which are subgraphs of $G'$. Dotted curves represent vertices and edges that existed in $G'$ but not in the cycle. Notice that every vertex and edge in $G'$ is included it at least one of these cycles. As before, blue and black curves represent edges in the graph while red curves represent collection of vertices and edges.}
    \label{fig:a'b'c'}
\end{figure}
Notice that $w,w_1,\ldots,w_k$ are in all of the cycles. Let $A,B,C$ be the graphs $A',B',C',$ respectively, after the edges $(w,w_i)$ for $1\leq i\leq k$ are added. Notice then that $A,B,C$ are starcle graphs (up to relabeling vertices appropriately) satisfying the conditions of \cref{starcle}. Define graph $D \in \{A,B,C\}.$
Let $D'$ be the star graph centered at $w$ with vertex set $V(D).$ We know that there exist $k$ disjoint paths between two vertices $\tau_1$ and $\tau_2$ in $\FS(D,D')$ if $\rho_1^{-1}(w) = \rho_2^{-1}(w) = w.$ 

Therefore, because $A,B,C$ share at least $k$ vertices aside from $w,$ namely the $k$ vertices on clockwise path from $z_1$ to $z_2$ or the clockwise path from $z_7$ to $z_8$ not including $w$, we know that by applying \cref{multifs} twice, there exist $k$ disjoint paths between $\tau_1$ and $\tau_2$ in $\FS(G^{\star\star\star},H^{\star\star})$ if $\rho_1^{-1}(w) = \rho_2^{-1}(w) = w.$

We now construct the $k$ disjoint paths from $\sigma$ to $\sigma \circ (m+1,m+2)$. We will start by creating an overarching structure. Let $\tau_1$ be a permutation on $[m+2]$ satisfying the following properties:
\begin{itemize}
    \item The vertices $\tau_1^{-1}(x_1), \tau_1^{-1}(x_2), \ldots, \tau_1^{-1}(x_\ell)$ appear consecutively in the large cycle of $G^\star$ such that $\tau_1^{-1}(x_1)=z_3$ and $\tau_1^{-1}(x_\ell)=z_5$,
    \item The vertices $\tau_1^{-1}(y_1), \tau_1^{-1}(y_2), \ldots, \tau_1^{-1}(y_\ell)$ appear consecutively in the large cycle of $G^\star$ such that $\tau_1^{-1}(y_1)=z_9$ and $\tau_1^{-1}(y_\ell)=z_{11}$,
    \item We have $\tau_1^{-1}(m+1)=m+1$ and $\tau_1^{-1}(m+2)=m+2$.
\end{itemize}
Define the following permutation as a composition of several transpositions.
\begin{align*}
\Sigma_1 = \, &(m+1,\tau_1^{-1}(x_1)) \circ (\tau_1^{-1}(x_1), \tau_1^{-1}(x_2) \circ \cdots \circ (\tau_1^{-1}(x_{\ell-1}) \circ \tau_1^{-1}(x_\ell)) \circ\\
&(m+2,\tau_1^{-1}(y_1))\circ (\tau_1^{-1}(y_1) \circ \tau_1^{-1}(y_2))\circ \cdots\circ (\tau_1^{-1}(y_{\ell-1},\tau_1^{-1}(y_\ell)).
\end{align*}
Let $\tau_2 = \tau_1 \circ \Sigma_1$ so that $\tau_2^{-1}(m+1)=z_5$ and $\tau_2^{-1}(m+2)=z_{11}.$

Let $\tau_3$ be a permutation where $\tau_3^{-1}(x_1) = m+2,$ $\tau_3^{-1}(y_1) = m+1,$ $\tau_3^{-1}(m+1)=z_5,$ and $\tau_3^{-1}(m+2)=z_{11}.$ Let $\tau_4 = \tau_3 \circ \Sigma_2$ where
$$\Sigma_2=(z_5, m+2) \circ (z_{11},m+1).$$
For all of these permutations $\tau \in \{\tau_1,\tau_2,\tau_3,\tau_4\},$ assume $\tau(w) = w.$ Define $\tau_{1,r}$ for $1\leq r \leq k$ to be distinct permutations of the form $\tau_1 \circ \pi$ where $\pi$ is a permutation of the vertices on the clockwise path from $z_1$ to $z_2.$ Similarly, define $\tau_{3,r}$ for $1\leq r\leq k$ to be distinct permutations of the form $\tau_3\circ \pi$ where $\pi$ is a permutation of the vertices on the clockwise path from $z_1$ to $z_2.$ Let $\tau_{2,r} = \tau_{1,r} \circ \Sigma_1,$ and let $\tau_{4,r} = \tau_{3,r}\circ \Sigma_2.$

Consider the paths from $\tau_{1,r}$ to $\tau_{2,r}$ by performing the transpositions comprising $\Sigma_1$ in order.
Notice that for distinct $r$ and $s$ the path from $\tau_{1,r}$ to $\tau_{2,r}$ and the path from $\tau_{1,s}$ to $\tau_{2,s}$ are disjoint as the positions of the people standing on vertices on the clockwise path from $z_1$ to $z_2$ are unchanged by the swaps of $\Sigma_1.$

Similarly the paths from $\tau_{3,r}$ to $\tau_{4,r}$ by performing the transpositions $\Sigma_2$ are also disjoint. Moreover, none of these paths can share a vertex with a path from $\tau_{1,s}$ to $\tau_{2,s}.$ In intermediate permutations of the path from $\tau_{3,r}$ to $\tau_{4,r},$ the position corresponding to $m+1$ is either $z_5$ or $m+2,$ and the position corresponding to $m+2$ is either $z_{11}$ or $m+1$ respectively.
This only happens in the path from $\tau_{1,s}$ to $\tau_{2,s}$ at $\tau_{2,s}$ where the position corresponding to $m+1$ is $z_5$ and the position corresponding to $m+2$ is $z_{11}.$ But notice that the position of $x_1$ is $m+1$ and the position of $y_1$ is $m+2$ in $\tau_{2,s}.$ This never happens in the path from $\tau_{3,r}$ to $\tau_{4,r}.$

Consider a vertex cut between $\sigma$ and $\sigma \circ (m+1,m+2)$ of $k-1$ vertices. Because the paths from $\tau_{1,r}$ to $\tau_{2,r}$ and the paths from $\tau_{3,r}$ to $\tau_{4,r}$ are all disjoint, if we remove $k-1$ vertices, there must always be some $r'$ for which there still exists a path between $\tau_{1,r'}$ to $\tau_{2,r'}$ and there exists a path between $\tau_{3,r'}$ and $\tau_{4,r'}.$

Note that in $\FS(G^{\star \star}, H^{\star\star})$ there are $k$ disjoint paths between $\sigma$ and $\tau_{1,r'}.$ These paths also exist in $\FS(G^\star,H^\star).$ Therefore removing $k-1$ vertices doesn't disconnect $\sigma$ and $\tau_{1,r'}.$ There must also exist $k$ disjoint paths between $\tau_{2,r'}$ and $\tau_{3,r'}$ in $\FS(G^{\star \star \star},H^{\star \star}),$ so again $\tau_{2,r'}$ and $\tau_{3,r'}$ must remain connected. Finally $\tau_{4,r'}$ and $\sigma\circ(m+1,m+2)$ must remain connected since there exist $k$ disjoint paths between them in $\FS(G^{\star \star}, H^{\star\star}).$
So we have a path between $\sigma$ and $\sigma\circ(m+1,m+2)$. Note that the claims of the existence of $k$ disjoint paths are due to \cref{starcle}.

Therefore, there exist $k$ disjoint paths between them by Menger's Theorem \cite{diestelblocks}.
\end{proof}
We now present the following lemma about the embeddability of $(G^\star,H^\star).$ Because the proof is almost identical to that of Wang and Chen \cite{wang2023connectivity}, we omit it.
\begin{lemma} [\protect{\cite[Lemma 3.5]{wang2023connectivity}}] 
    Let $n$ be a large positive integer, and let $m,\ell,G^{\star\star},$ and $H^{\star\star}$ be as described above. Let $\Gamma = \{x_1,\ldots,x_\ell,y_1,\ldots,y_\ell,z_3,z_5,z_9,z_{11}\}.$ Let $q_i = \lfloor p_0n/(5\ell) \rfloor$ for all $i \in \Gamma,$ and let $q_i = \lfloor n/(2m) \rfloor$ for all $i \in [m]\setminus \Gamma$ where $p_0$ be a probability such that
    $$p_0\geq \frac{\exp((k+7)/4\cdot\log n)^{2/3}}{n^{1/2}}.$$
    Let $p_1$ and $p_2$ be probabilities such that $p_1p_2 \geq p_0^2$ and $p_2\geq p_1 \geq p_0/\ell.$
    Let $X$ and $Y$ be independently-chosen random graphs in $\mathcal{G}(n,p_1)$ and $\mathcal{G}(n,p_2),$ respectively. If $n$ is sufficiently large, then the probability that the pair $(G^{\star\star},H^{\star\star})$ is $(q_1,\ldots,q_m)$-embeddable in $(X,Y)$ is at least $1-n^{-n/3}.$
\end{lemma}
As a result of these lemmas, we have \cref{wangchenkconnectedwhp} which completes the proof of \cref{kconnectivitywhp}. We omit its proof due to its similarity to the proof in Wang and Chen \cite{wang2023connectivity}.
\section{Future Research Directions}
Relating to \cref{genwilson}, we think this result can be generalized to the following two conjectures. There is already some evidence of this with \cref{starcle}, which discusses starcle graphs, as we found that the number of disjoint paths between $\sigma$ and $\rho$ from the lemma statement was significantly higher than the minimum degree of the starcle graph.
\begin{conjecture}\label{starvertexconj}
    Let $G = \FS(X,\Star_n)$ where $X$ is a graph on $n$ vertices. Assume $G$ is connected. If $\sigma$ and $\rho$ are permutations in $V(G),$ the maximum number of disjoint paths between $\sigma$ and $\rho$ is $\min\{\deg(\sigma),\deg(\rho)\}.$ 
\end{conjecture}

Another interesting direction is to consider the following question posed by Colin Defant.
\begin{question}
    What are the connectivities of the connected components of $\FS(X,P_n)$ or $\FS(X,C_n)$ where $X$ is a graph on $n$ vertices?
\end{question}

For these graphs, the method of atomic parts does not seem to work as there are very few useful graph automorphisms. 

In \cref{kban1,kban2,kconnectivitywhp}, we fixed $k$ and varied $n$ and found that similar bounding conditions to the analogs in classical connectivity show $\FS(X,Y)$ is $k$-connected. The following question considers letting $k$ scale with $n.$
\begin{problem}
    Choose $k = k(n)$ to be some function on $n$. What are the bounds on the minimum degree to ensure $\FS(X,Y)$ is $k$-connected as in \cref{kban1,kban2}. What are the bounds on the threshold probability as in \cref{kconnectivitywhp}?
\end{problem}

\section*{Acknowledgements}
We would like to thank Dr.\ Tanya Khovanova of the Department of Mathematics at the Massachusetts Institute of Technology (MIT) for her advice, and Dr.\ Colin Defant of the Department of Mathematics at Harvard University for useful discussions. We thank Dr.\ Slava Gerovitch and Professor Pavel Etingof of the Department of Mathematics at MIT for creating and operating the PRIMES-USA program which has given this opportunity. We also thank the MIT PRIMES-USA Program as well as the Department of Mathematics at MIT.
\bibliographystyle{amsinit}
\bibliography{ref}

\providecommand{\bysame}{\leavevmode\hbox to3em{\hrulefill}\thinspace}
\providecommand{\MR}{\relax\ifhmode\unskip\space\fi MR }
\providecommand{\MRhref}[2]{%
  \href{http://www.ams.org/mathscinet-getitem?mr=#1}{#2}
}
\providecommand{\href}[2]{#2}
\begin{thebibliography}{10}

\bibitem{alon2021extremal}
N. Alon, C. Defant, and N. Kravitz, \emph{Typical and extremal aspects of friends-and-strangers graphs}, J. Combin. Theory Ser. B \textbf{158} (2023), 3--42. \MR{4513816} \doi{10.1016/j.jctb.2022.03.001}

\bibitem{bangachev2022asymmetric}
K. Bangachev, \emph{On the asymmetric generalizations of two extremal questions on friends-and-strangers graphs}, European J. Combin. \textbf{104} (2022), Paper No. 103529, 26. \MR{4400016} \doi{10.1016/j.ejc.2022.103529}

\bibitem{bollobas2016packing}
B. Bollob\'as, S. Janson, and A. Scott, \emph{Packing random graphs and hypergraphs}, Random Structures Algorithms \textbf{51} (2017), no.~1, 3--13. \MR{3668844} \doi{10.1002/rsa.20673}

\bibitem{defant2021friends}
C. Defant and N. Kravitz, \emph{Friends and strangers walking on graphs}, Comb. Theory \textbf{1} (2021), Paper No. 6, 34. \MR{4396211} \doi{10.5070/C61055363}

\bibitem{diestelblocks}
R. Diestel, \emph{Graph theory}, fifth ed., Graduate Texts in Mathematics, vol. 173, Springer, Berlin, 2018. \MR{3822066}

\bibitem{godsil1981cayley}
C.~D. Godsil, \emph{Connectivity of minimal {C}ayley graphs}, Arch. Math. (Basel) \textbf{37} (1981), no.~5, 473--476. \MR{643291} \doi{10.1007/BF01234384}

\bibitem{goring2000digraphmenger}
F. G\"oring, \emph{Short proof of {M}enger's theorem}, Discrete Math. \textbf{219} (2000), no.~1-3, 295--296. \MR{1761733} \doi{10.1016/S0012-365X(00)00088-1}

\bibitem{jeong2022diameters}
R. Jeong, \emph{On the diameters of friends-and-strangers graphs}, Comb. Theory \textbf{4} (2024), no.~2, Paper No. 2. \doi{10.5070/C64264229}

\bibitem{menger1927theorem}
K. Menger, \emph{Zur allgemeinen kurventheorie}, Fund. Math. \textbf{10} (1927), no.~1, 96--115.

\bibitem{milojevic2022connectivity}
A. Milojevi\'c, \emph{Connectivity of old and new models of friends-and-strangers graphs}, Adv. in Appl. Math. \textbf{155} (2024), Paper No. 102668, 53. \MR{4689232} \doi{10.1016/j.aam.2023.102668}

\bibitem{wang2023connectivity}
L. Wang and Y. Chen, \emph{Connectivity of friends-and-strangers graphs on random pairs}, Discrete Math. \textbf{346} (2023), no.~3, Paper No. 113266, 10. \MR{4513695} \doi{10.1016/j.disc.2022.113266}

\bibitem{watkins1970atomicparts}
M.~E. Watkins, \emph{Connectivity of transitive graphs}, J. Combinatorial Theory \textbf{8} (1970), 23--29. \MR{266804}

\bibitem{wilson1974graph}
R.~M. Wilson, \emph{Graph puzzles, homotopy, and the alternating group}, J. Combinatorial Theory Ser. B \textbf{16} (1974), 86--96. \MR{332555} \doi{10.1016/0095-8956(74)90098-7}

\end{thebibliography}

\end{document}